\RequirePackage{fix-cm}
\documentclass[twocolumn]{svjour3}          
\smartqed  
\usepackage{graphicx}
\usepackage{amssymb}
\usepackage{amscd}
\usepackage{amsmath}
\DeclareMathAlphabet\gothic{U}{euf}{m}{n}

\usepackage{algorithm2e}

\usepackage{todonotes}

\usepackage{booktabs}

\DeclareMathAlphabet\gothic{U}{euf}{m}{n}

\begin{document}

\title{Nilpotent Approximations of Sub-Riemannian Distances for Fast Perceptual Grouping of Blood Vessels in 2D and 3D
}


\author{Erik J. Bekkers         \and
        Da Chen                 \and
        Jorg M. Portegies
}


\institute{E.J. Bekkers \and J.M. Portegies \at
              Centre for Analysis, Scientific computing and Applications (CASA),
              Eindhoven University of Technology, the Netherlands \\
              \email{e.j.bekkers@tue.nl, j.m.portegies@tue.nl}\\
           \and
           D. Chen \at
              University Paris Dauphine, PSL Research University\\
              CNRS, UMR 7534, CEREMADE, 75016 Paris, France \\
              \email{chenda@ceremade.dauphine.fr}
           }

\date{Received: date / Accepted: date}

\maketitle

\begin{abstract}
We propose an efficient approach for the grouping of local orientations (points on vessels) via nilpotent approximations of sub-Riemannian distances in the 2D and 3D roto-translation groups $SE(2)$ and $SE(3)$.
In our distance approximations we consider homogeneous norms on nilpotent groups that locally approximate $SE(n)$, and which are obtained via the exponential and logarithmic map on $SE(n)$. In a qualitative validation we show that the norms provide accurate approximations of the true sub-Riemannian distances, and we discuss their relations to the fundamental solution of the sub-Laplacian on $SE(n)$. The quantitative experiments further confirm the accuracy of the approximations. Quantitative results are obtained by evaluating perceptual grouping performance of retinal blood vessels in 2D images and curves in challenging 3D synthetic volumes. The results show that 1) sub-Riemannian geometry is essential in achieving top performance and 2) that grouping via the fast analytic approximations performs almost equally, or better, than data-adaptive fast marching approaches on $\mathbb{R}^n$ and $SE(n)$.
\end{abstract}

\section{Introduction}
\label{intro}
\begin{figure}
  \centering
  \includegraphics[width=0.48\textwidth]{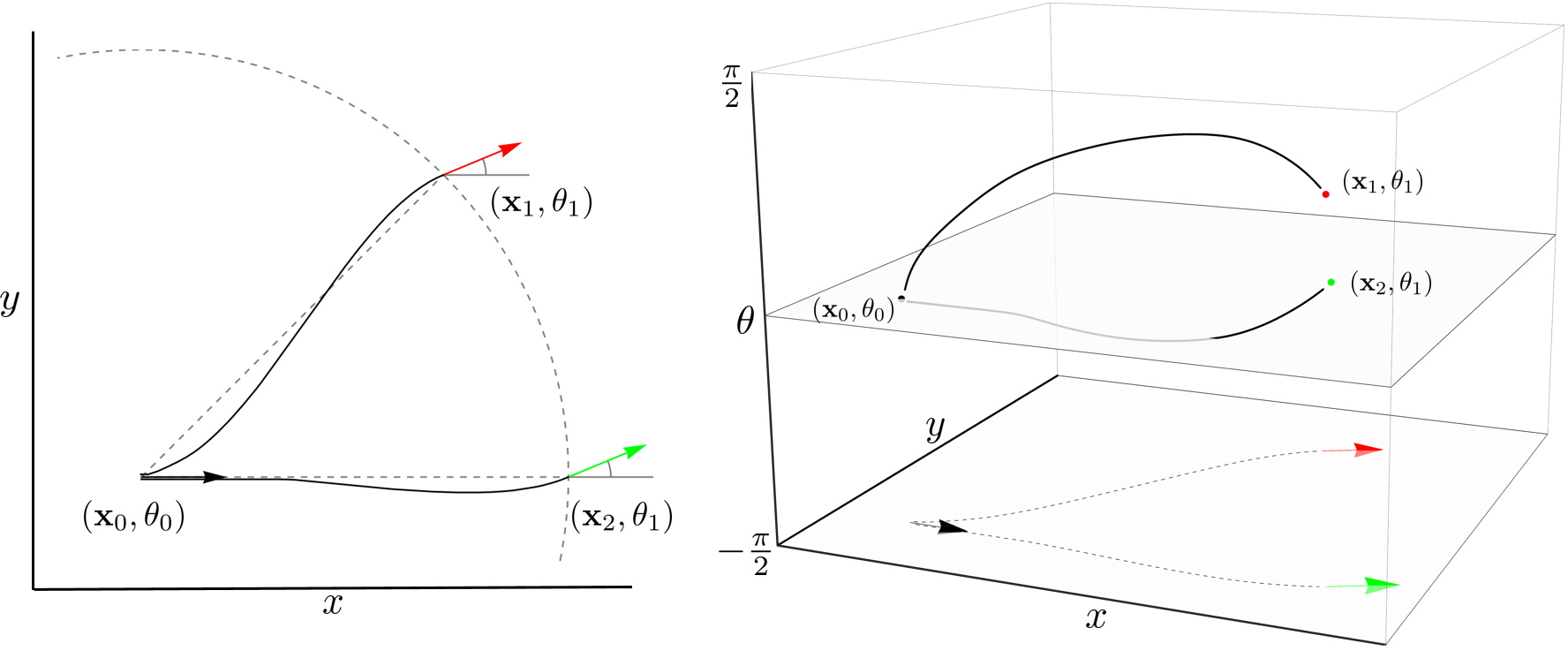}
  \caption{
   The red and green arrow have equal spatial and angular distance to the origin (black arrow). In a flat geometry on $\mathbb{R}^2 \times S^1$ the distance between the red and green arrow and the source would be equal, and the geodesics straight lines (see dashed lines). In sub-Riemannian geometry on $SE(2)$ the green arrow has a shorter distance to the source. The left image shows 2D projections of the sub-Riemannian geodesics in solid black, the right image shows their paths in $SE(2)$.}\label{fig:AlignmentSE2}
\end{figure}

In this paper we derive analytic formulas for approximations of sub-Riemannian distances on the 2D and 3D rotation translation groups, denoted respectively with $SE(2)$ and $SE(3)$. Additionally, we extend the perceptual grouping algorithm \cite{Cohen2001multiple} for clustering of local orientations (points on blood vessels). Here clustering is based on alignment of local orientations, which is quantified using sub-Riemannian distances on $SE(n)$, see Fig.~\ref{fig:AlignmentSE2} for an illustration.

\subsection{Nilpotent Approximation}
The sub-Riemannian distances on $SE(n)$ are approximated via norms on the vectors obtained from the logarithmic map (from group elements to the Lie algebra).
This approach is motivated by problems from sub-Riemannian geometry in nilpotent Lie groups, in which such homogenous norms provide exact fundamental solutions to sub-Laplacians.

The vectors obtained by the logarithmic map, expressed in a left-invariant basis, are the so-called exponential coordinates of the first kind. For a nilpotent group of step two, like the Heisenberg group, these coordinates define (together with a group product defined via the Baker-Campbell-Hausdorf (BCH) formula) a \emph{global} isomorphism to the group. In our $SE(n)$ setting we have to truncate the commutator series in the BCH formula due to non-vanishing (higher-order) commutators, yielding a corresponding Heisenberg type approximation which we denote with $(SE(n))_0$.
The obtained Taylor development of the group product and associated left-invariant vector fields gives rise to a \emph{local} approximation of the (sub-Riemannian) flows on $SE(2)$ in the sense of Rothschild and Stein \cite{rothschild-stein}.

We then define a norm on $(SE(n))_0$ based on the Folland-Kaplan-Kor{\'a}nyi gauge, which is known for its relation to the fundamental solution of the sub-Laplacian on the Heisenberg group \cite{folland1973fundamental,kaplan1980fundamental,koranyi1982kelvin}.
We reason that the Folland-Kaplan-Kor{\'a}nyi provides an accurate approximation to the fundamental solution on $SE(n)$  as well, as it provides the exact fundamental solution on the Heisenberg type approximation $(SE(n))_0$. As such, we provide an approach to approximating the heat kernel and fundamental solution of the sub-Laplacian on $SE(n)$, as an alternative to the works \cite{Duits2010,portegies2016arxiv,citti2006cortical}.

The distance associated with the Folland-Kaplan-Kor{\'a}nyi type norm on $(SE(n))_0$ is locally equivalent to the sub-Riemannian distance on $SE(2)$, as was formally proven in full generality in the seminal work by Nagel, Stein and Wainger \cite{nagel-stein-wainger}. In this paper, we show by qualitative and quantitative comparison that the norm on $(SE(n))_0$ indeed provides a sharp \emph{local} approximation of the sub-Riemannian distances on $SE(n)$.


\subsection{Perceptual Grouping}

The motivation for perceptual grouping of local orientations comes from problems in medical image analysis in which the topologically correct reconstruction of vessel (and pulmonary) trees is of great importance in biomarker research and surgery planning.
Knowing the correct connectivity in tree structures not only allows for local biomarker analysis (e.g., studies on bifurcation and crossing properties \cite{Leontidis2016}), but also allows for higher level biomarker research via statistics on tree structures \cite{feragen2015geodesic}.
Topological knowledge of vessel trees is also essential in determining artery/vein classification problems \cite{estrada2015retinal,Eppenhof2015,Dashtbozorg2014}. 
Finally, in many medical applications involving vessel analysis, including topological tree reconstruction, distances between local orientations play a crucial role \cite{de2016graph,turetken2016reconstructing,lo2010vessel,shang2011vascular,Abbasi-Sureshjani2017,favali2016analysis}. 
The approximate sub-Riemannian distance in this paper is analytic, fast, and easy to implement, and as such may be a useful tool for algorithms that rely on local orientation analysis.

Sub-Riemannian models are shown to be effective in both image processing and in neuropsychological models for line perception in the primary visual cortex \cite{Petitot2003,citti2006cortical,Sarti2015,yuriSE2FINAL,mashtakov2016cortical,Duits2013,boscain2014,Bekkers2015SIIMS,prandi2015highly,favali2016analysis}. In this paper we indeed observe by quantitative validation of automatic connectivity analysis that sub-Riemannian distances are preferred over their (full) Riemannian counter parts.

\begin{figure*}
  \centering
  \includegraphics[width=1\textwidth]{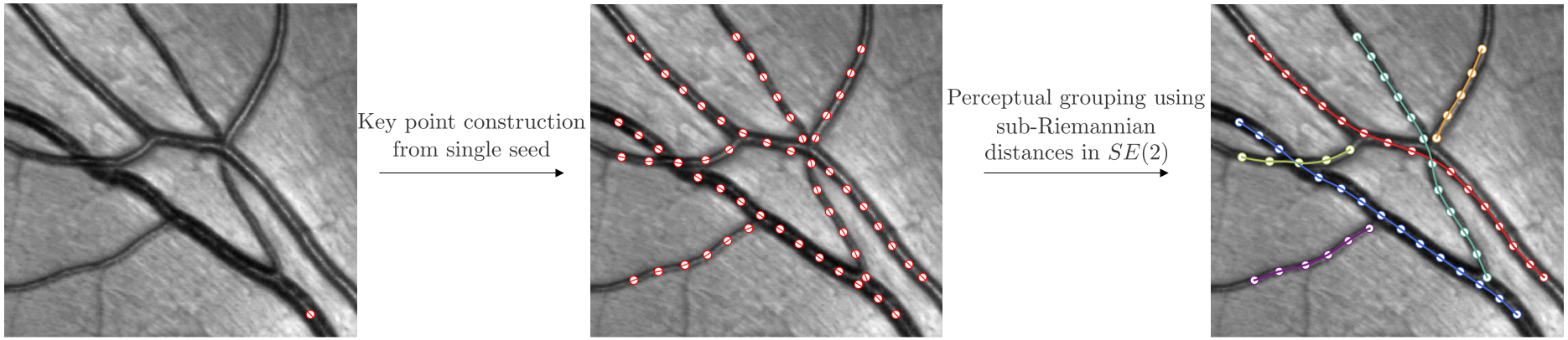}
  \caption{The pipeline for grouping vessel segments consists of 2 steps. First, key points are generated (from a single source point) using minimal path tracking with key points \cite{benmansour2009fast}. Second, the automatically generated key points, with estimated orientations, are grouped based on an adaption of the perceptual grouping algorithm \cite{Cohen2001multiple} with the use of sub-Riemannian distances on $SE(2)$. The result on the right is obtained with the nilpotent approximations of the sub-Riemannian distances in $SE(2)$. }\label{fig:PipelineSE2}
\end{figure*}

The approach taken in this paper for doing connectivity analysis is based on the perceptual grouping algorithm proposed by Cohen \cite{Cohen2001multiple}. This algorithm turns a set of key points into a graph by iteratively adding edges between nodes based on their geodesic distances while putting constraints on the number of connections per node. The input set of key points may be obtained via key point tracking algorithms \cite{benmansour2009fast,kaul2012detecting,chen2016vessel}, as is done also in this paper, see Fig.~\ref{fig:PipelineSE2}.

In \cite{Cohen2001multiple} an isotropic metric was used to define the geodesic distances. Later, the perceptual grouping algorithm was adapted for use with anisotropic Riemannian metrics by Bougleux et al. \cite{Bougleux2008anisotropic}. In recent work \cite{Chen2017}, it was further extended for the grouping of $n$ closed contours for an a-priori $n$. There, a (sub-)Finsler metric on position orientation space was used, similar to the sub-Riemannian metric used in this paper. As in \cite{Bougleux2008anisotropic} and \cite{Chen2017}, we use the main algorithm of \cite{Cohen2001multiple} as a backbone, but we change the metric used for perceptual grouping and we impose an additional constraint to avoid closed loops (which are physically not realistic in the vessel networks of interest).

With quantitative experiments we show that perceptual grouping with sub-Riemannian distances on $SE(n)$ is preferred over the use of (full) Riemannian distances on $SE(n)$, which is in turn preferred over grouping with distances on $\mathbb{R}^n$. Furthermore, the analytic approximations allow for fast perceptual grouping with competitive performance compared to \emph{data-adaptive} sub-Riemannian distances computed via fast marching.

\subsection{Paper Outline}
In Sec.~\ref{sec:SE2} and Sec.~\ref{sec:SE3} we derive approximations for sub-Riemannian distances in respectively $SE(2)$ and $SE(3)$. There, for each Lie group we first provide the preliminaries, then define the sub-Riemannian distance, and then describe the proposed approximations. In Sec.~\ref{sec:algorithms} the algorithms (perceptual grouping, fast marching and key point tracking) are described, including an overview of the different distances used in this paper. In Sec.~\ref{sec:experiments} we then compare the performance of the perceptual grouping algorithm using different distances, first on $\mathbb{R}^2$ and $SE(2)$ in Subsec.~\ref{subsec:PGSE2}, then on $\mathbb{R}^3$ and $SE(3)$ in Subsec.~\ref{subsec:PGSE3}. General conclusions are provided in Sec.~\ref{sec:conclusion}.

\section{Sub-Riemannian Distance and its Approximation in $SE(2)$}\label{sec:SE2}
\subsection{The Lie Group $SE(2)$}
\subsubsection{SE(2)}\label{subsec:SE2}
In order to measure distances between local orientations we will consider the Lie group SE(2) as our base manifold.
The group $SE(2) = \mathbb{R}^2 \rtimes SO(2)$ is the semi-direct product of the group of planar translations $\mathbb{R}^2$ and rotations $SO(2)$, and its group product and inverse are respectively defined via:
\begin{equation}\label{eq:groupproduct}
\begin{array}{rl}
g \cdot g' &= ( \mathbf{x},\mathbf{R}_\theta ) \cdot ( \mathbf{x}', \mathbf{R}_{\theta'} ) 
= ( \mathbf{R}_\theta \mathbf{x}' + \mathbf{x}, \mathbf{R}_{\theta + \theta'}),\\
g^{-1} &= ( -\mathbf{R}_{\theta}^{-1} \mathbf{x} , \mathbf{R}_{\theta}^{-1}),
\end{array}
\end{equation}
with group elements $g,g' \in SE(2)$.
The group acts on the (coupled) space of positions and orientations $\mathbb{R}^2 \rtimes S^1$ via
$$
g \cdot (\mathbf{x}',\theta') = ( \mathbf{R}_\theta \mathbf{x}' + \mathbf{x}, \theta +\theta').
$$
Since $( \mathbf{x} , \mathbf{R}_\theta ) \cdot ( \mathbf{0}, 0 ) = ( \mathbf{x} , \theta )$, we can uniquely identify the roto-translation group $SE(2)$ with the space of positions and orientations $\mathbb{R}^2 \rtimes S^1$.

\subsubsection{The Lie Algebra, Exponential Map and Commutators}
%
The Lie algebra associated with $SE(2)$ is the real vector space $\gothic{se}(2) = \operatorname{span} \{ A_1, A_2, A_3 \}$ together with a bilinear operator $[\cdot,\cdot]:\gothic{se}(2)\times\gothic{se}(2)\rightarrow\gothic{se}(2)$ called the Lie bracket (which we define below in Eq.~(\ref{eq:LieBracket})).
The generators of the Lie algebra are given by the differential frame $\left.\{\partial_\theta,\partial_x,\partial_y\}\right|_{(0,0,0)}$ at the origin
\begin{equation}\label{eq:LieAlgebra}
A_1 = \left.\partial_\theta\right|_{(0,0,0)}, \;\; A_2 = \left.\partial_x\right|_{(0,0,0)}, \;\; A_3 = \left.\partial_y\right|_{(0,0,0)},
\end{equation}
which define corresponding left-invariant vector fields
\begin{equation}\label{eq:LeftInvVectorFields}
\begin{split}
  \left.\mathcal{A}_1\right|_{g} = (L_{g})_* A_1 &= \hspace{0.8em}\left.\partial_\theta\right|_{g}, \\
  \left.\mathcal{A}_2\right|_{g} = (L_{g})_* A_2 &= \hspace{0.8em}\cos \theta \left.\partial_x\right|_{g} + \sin \theta \left.\partial_y\right|_{g}, \\
  \left.\mathcal{A}_3\right|_{g} = (L_{g})_* A_3 &= -\sin \theta \left.\partial_x\right|_{g} + \cos \theta \left.\partial_y\right|_{g}
\end{split}
\end{equation}
via the push-forward of left-multiplication, denoted by $(L_g)_*$, and with $g = (x,y,\theta) \in SE(2)$.

The exponential map $\operatorname{Exp}: \gothic{se}(2) \rightarrow SE(2)$ defines a mapping from a vector $X \in \gothic{se}(2)$ in the tangent space at $g=(0,0,0)$ to an element in the group $SE(2)$ by following an integral curve along the left-invariant vector field $(L_g)_* X$. 
The logarithmic map $\operatorname{Log}: SE(2) \rightarrow \gothic{se}(2)$ defines the mapping from group element to tangent vector at $g=(0,0,0)$.

The Lie bracket for vector fields is defined as follows
\begin{equation}\label{eq:LieBracket}
\begin{split}
  [X,Y]:&= \underset{t\rightarrow 0}{\operatorname{lim}} \frac{\gamma(t) - e}{t^2}, \;\;\;\;\;\;\; \text{with} \\
  \gamma(t) &= \operatorname{Exp}(-t Y)\operatorname{Exp}(-t X)\operatorname{Exp}(t Y)\operatorname{Exp}(t X).
\end{split}
\end{equation}
I.e., it describes the infinitesimal displacement by following a path moving forth and back in $X$ and $Y$ directions. The Lie bracket of two vectors defines a new vector (the commutator) and the Lie bracket of two vector fields defines a new vector field. The non-zero commutators of $\gothic{se}(2)$ are
\begin{equation}\label{eq:commutatorsSE2}
\begin{split}
    [A_1,A_2] &= -[A_2,A_1] = \;\;\;A_3, \\
    [A_1,A_3] &= -[A_3,A_1] = -A_2.
\end{split}
\end{equation}


\subsection{Sub-Riemannian Geometry in SE(2)}\label{subsec:SRGeometrySE2}
We consider a sub-Riemannian geometry on $SE(2)$ by measuring distances between two points in $SE(2)$ via the lengths of shortest horizontal paths. A horizontal path is a curve $\gamma:[t_0,t_1]\subset\mathbb{R}\rightarrow SE(2)$ with tangent vectors $\dot{\gamma}(t) \in \left.\Delta\right|_{\gamma(t)}:= \operatorname{span} \{ \left. \mathcal{A}_1\right|_{\gamma(t)}, \left. \mathcal{A}_2 \right|_{\gamma(t)}\}$, where $\Delta$ denotes the sub-Bundle of the full tangent bundle $T(SE(2)):= \operatorname{span} \{\mathcal{A}_1,\mathcal{A}_2,\mathcal{A}_3\}$. Lengths of horizontal curves with $\dot{\gamma}(t) = u^1(t) \left.\mathcal{A}_1\right|_{\gamma(t)} + u^2(t) \left.\mathcal{A}_2\right|_{\gamma(t)}$ are measured by the sub-Riemannian metric tensor\footnote{Due to the fact the the metric tensor is degenerate in the $\mathcal{A}_3$ direction (tangent vectors are always contained within $\Delta$) it is not possible to represent the metric tensor in a standard form as an \emph{invertible} symmetric $3\times3$ matrix. This is however possible when including an additional term $\epsilon^{-2} \xi^2 |u^3|^2(t)$ after which the tensor becomes (anisotropic) Riemannian \cite{citti2006cortical,Sanguinetti2015CIARP}. This Riemannian approximation converges to the sub-Riemannian tensor when $\epsilon\rightarrow 0$ \cite[App.~A]{ThesisChenDa2016} and \cite[Thm.~2]{duits2016optimal}.}
\begin{equation}\label{eq:SRMetricTensorSE2}
  \left.\mathcal{G}^{\xi,C}\right|_{\gamma(t)} (\dot{\gamma}(t),\dot{\gamma}(t)):= C(\gamma(t))^2 (|u^1(t)|^2 + \xi |u^2(t)|^2),
\end{equation}
in which $C:SE(2)\rightarrow \mathbb{R}^+$ is an external cost which penalizes the curves to move through certain regions in $SE(2)$, $\xi$ is a parameter which balances the penalty of motion in the angular and spatial directions and has dimensions [1/length], and $u^1$ and $u^2$ are the control parameters of the curve $\gamma$.

The sub-Riemannian distances between two points $g_1,g_2 \in SE(2)$ is then given by
\begin{equation}\label{eq:SRDistancesSE2}
  d_0(g_1,g_2) := \operatorname{inf} \left\{ \int_0^1 \sqrt{ \left.\mathcal{G}^{\xi,C}\right|_{\gamma(t)}(\dot{\gamma}(t),\dot{\gamma}(t)) } {\rm d}t \right\},
\end{equation}
where the infimum is taken over Lipschitz continuous curves $\gamma \in \operatorname{Lip}([0,T],SE(2))$ with $\gamma(0) = g_1$, $\gamma(1) = g_2$, and $\dot{\gamma}(t) = u^1(t) \left.\mathcal{A}_1\right|_{\gamma(t)} + u^2(t) \left.\mathcal{A}_2\right|_{\gamma(t)}$. Note that due to the inclusion of an external cost function $C$ the distance $d$ is not strictly left-invariant, however, when substituting $C$ by $C_g:=C(g^{-1} h)$ in (\ref{eq:SRDistancesSE2}) we do have left-invariance (i.e., then $d(g \cdot g_1, g \cdot g_2) = d(g_1, g_2)$).

\subsection{A Nilpotent Approximation $(SE(2))_0$ of $SE(2)$}\label{sec:ApproximationSE2}

\subsubsection{A Local Approximation via the Baker-Campbell-Hausdorff Formula}
Consider the exponential map from Lie algebra $\gothic{se}(2)$ to the group $SE(2)$
\begin{equation}\label{eq:coordsfirst}
(c^1,c^2,c^3)\mapsto (x,y,\theta) = \operatorname{Exp}( c^1 A_1 + c^2 A_2 + c^3 A_3 ),
\end{equation}
with $\{A_i\}_{i=1}^3$ the basis vectors of $\gothic{se}(2)$ given in (\ref{eq:LieAlgebra}), and with $(c^1,c^2,c^3)$ the canonical coordinates of the first kind given by
%
%
\begin{equation} \label{eq:coordinatesfirstkind}
\begin{array}{l}
c^1 = \theta, \;\;\;\;\; c^2 = \left\{ \begin{array}{ll}
\tfrac{1}{2} \theta (y + x \cot \tfrac{\theta}{2} ) & \;\;\;\;\;\; \text{if} \;\;\; \theta \neq 0\\
x & \;\;\;\;\;\; \text{if} \;\;\; \theta = 0\\
\end{array}\right. ,\\
\;\;\;\;\;\;\;\;\;\;\;\;\;\;\;\;\;c^3 = \left\{ \begin{array}{ll}
\tfrac{1}{2} \theta (-x + y \cot \tfrac{\theta}{2} ) & \;\;\; \text{if} \;\;\; \theta \neq 0\\
y & \;\;\; \text{if} \;\;\; \theta = 0\\
\end{array}\right. .
\end{array}
\end{equation}
For two left-invariant vector fields $X = \sum_{i=1}^3 x^i \mathcal{A}_i$ and $Y = \sum_{i=1}^3 y^i \mathcal{A}_i$
the Baker-Campbell-Hausdorff (BCH) formula (see e.g. \cite{rossmann2002lie}) gives:
\begin{equation}\label{eq:BCH}
\begin{array}{rl}
\operatorname{Log}( \operatorname{Exp}(X) \operatorname{Exp}(Y) ) = &X + Y + \frac{1}{2} [X,Y] \\
&+ \frac{1}{12}([X,[X,Y]] + [Y,[Y,X]])\\
&\;\;\;\;\;\;\;\;\;\;\;\;\;\;\;\;\;\;\;\;\;\; + \mathcal{O}([\cdot,[\cdot,[\cdot,\cdot]]]),
\end{array}
\end{equation}
where $\mathcal{O}([\cdot,[\cdot,[\cdot,\cdot]]])$ denotes higher order nested brackets. Since the Lie algebra $\gothic{se}(2)$ is not nilpotent it has non-vanishing Lie brackets of order $\ge 2$ (cf. the commutator relations in (\ref{eq:commutatorsSE2})) the BCH formula gives an infinite series of nested Lie brackets.

\begin{figure*}
  \centering
  \includegraphics[width=\textwidth]{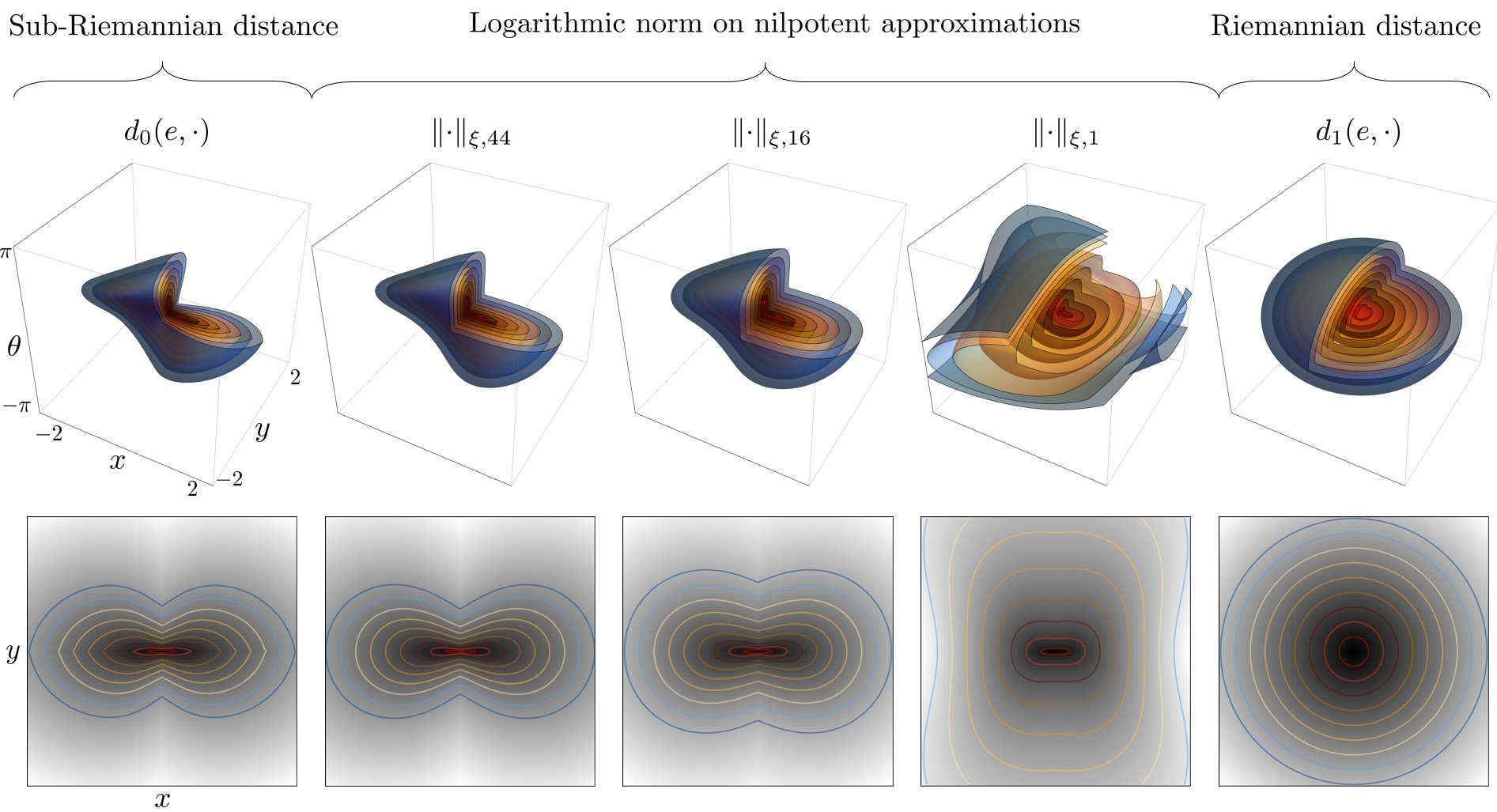}
  \caption{Distances on $SE(2)$ for $\xi=1$, $C=1$. Top row: Level sets of the distance volumes on $SE(2)$. Bottom row: Minimum intensity projections of the distances to the plane $\mathbb{R}^2$ with level set contours. From left to right: The sub-Riemannian distance $d_{0}(e,\cdot)$, see Eq.~(\ref{eq:SRDistancesSE2}); Homogenous norms $\lVert \cdot \rVert_{\xi,\zeta}$, see Eq.~(\ref{eq:approxSE2}), of the nilpotent approximation $(SE(2))_0$ for respectively $\zeta=44$, $\zeta=16$ (Folland-Kaplan-Kor{\'a}nyi gauge) and $c=\zeta$; The ($\xi$-isotropic) Riemannian distance $d_{1}(e,\cdot)$ on $SE(2)$, see Table~\ref{tab:metrics} for an overview of the different distances. }\label{fig:DistancesSE2}
\end{figure*}

Here, we approximate the BCH formula $SE(2)$ as\footnote{Note that such approximations of the BCH formula were already introduced in \cite[Thm.~2.22]{nagel-stein-wainger} in the general setting by Nagel, Stein and Wainger in 1985.}
\begin{equation}
\operatorname{Log}( \operatorname{Exp}(X) \operatorname{Exp}(Y) ) \approx X + Y + \frac{1}{2} [X,Y],
\end{equation}
by omitting the Lie brackets of order 2 (once nested brackets) and higher, as if our Lie algebra $\gothic{se}(2)$ is nilpotent of step 2.
Then, together with the commutator relations $[A_i,A_i]=0$, $A_3 = [A_1,A_2]$, and again omitting Lie brackets of order 2 (i.e., setting $[A_1,A_3] = [A_1,[A_1,A_2]] = 0$), the BCH formula defines a group product on the vector space $\mathbb{R}^3$ of the canonical coordinates of the first kind via
\begin{multline} \label{eq:GroupLawNilpotent}
 (x^1,x^2,x^3) \cdot (y^1,y^2,y^3) \\= \left(x^1+y^1,x^2+y^2,x^3+y^3 + \frac{1}{2}(x^1y^2-x^2y^1)\right).
\end{multline}
The new group product (\ref{eq:GroupLawNilpotent}), where the elements are expressed in coordinates of the first kind (cf. Eq.~(\ref{eq:coordsfirst})),
gives rise to a nilpotent Heisenberg group. It is a local\footnote{With $g_1,g_2 \in SE(2)$ chosen close enough such that higher order terms in (\ref{eq:BCH}) can be neglected.} approximation of the true group product $g_1 \cdot g_2 = \operatorname{Exp}(\sum_{i=1}^3 x^i A_i) \cdot \operatorname{Exp}(\sum_{i=1}^3 y^i A_i)$ given by (\ref{eq:groupproduct}).
We denote this group by $(SE(2))_0 = H(3)$, with $H(3)$ the 3 dimensional (nilpotent) Heisenberg group. Note that if $(x^1,x^2,x^3)$ and $(y^1,y^2,y^3)$ were coordinates of the first kind for a group with a step-2 nilpotent algebra, then this new group would be globally isomorphic to that group. The new group $(SE(2))_{0}$ defines a homogeneous Carnot group with respect to the dilations
\begin{equation}\label{eq:Dilations}
  \delta_s(\mathbf{c}) = (s \, c^1, s \, c^2, s^2 \, c^3).
\end{equation}



\subsubsection{Homogeneous Norms on $(SE(2))_0$ and the Fundamental Solution of the sub-Laplacian}
In our approximation of the sub-Riemannian distance $d_{0}$ of Eq.~(\ref{eq:SRDistancesSE2}) we make use of the following homogenous norm on $(SE(2))_{0}$:
\begin{equation}\label{eq:FKKGauge}
\lVert \mathbf{c} \rVert_\zeta := \sqrt[4]{(|c^1|^2 + |c^2|^2)^2 + \zeta \, |c^3|^2},
\end{equation}
with constant $\zeta>0$ a relative penalty for the non-horizontal part.
For $\zeta=16$ this norm coincides with the well-known Folland-Kaplan-Kor{\'a}nyi gauge, which is a widely studied norm on Carnot groups due to its relation to fundamental solutions of sub-Laplacians \cite{bonfiglioli2007stratified}:
Folland found that $\lVert \textbf{c} \rVert_{16}^{2-Q}$, with homogeneous dimensions $Q$, is (a constant multiple of) the fundamental solution of the canonical sub-Laplacian
on the Heisenberg group \cite{folland1973fundamental}; Kaplan showed that this relation more generally holds for H-type (Carnot) groups \cite{kaplan1980fundamental}; Kor{\'a}nyi derived many more of its properties in relation to harmonic analysis and potential theory \cite{koranyi1982kelvin}.

In relation to sub-Riemannian geometry on $SE(2)$ and its sub-Laplacian $\mathcal{L}:=\mathcal{A}_1^2 + \mathcal{A}_2^2$, we find that the fundamental solution $\Gamma$ of $\mathcal{L}$ can be approximated by the (explicit) fundamental solution of the canonical sub-Laplacian $\mathcal{L}_0: = \mathcal{X}_1^2 + \mathcal{X}_2^2$, with Jacobian basis
$
\mathcal{X}_1 = \partial_{c^1} + \frac{c^2}{2}\partial_{c^3}$, $\mathcal{X}_2 = \partial_{c^2} - \frac{c^1}{2}\partial_{c^3}
$ on $(SE(2))_{0}$. This solution in fact coincides with one of the approximations of $\Gamma$
found by Duits \& Franken \cite{Duits2010}. There, the fundamental solution of $\mathcal{L}$ was first approximated by relying on a contraction of $SE(2)$ to a 3-dimensional Heisenberg group (via dilations on the group $SE(2)$), and then derived the Gaussian estimates based on the homogeneous norm $\lVert \cdot \rVert_1$, i.e., $\zeta=1$, with exponential coordinates derived from the contraction.

In our study on the sub-Riemannian distance approximations we found that even sharper estimates could be obtained by relying on the explicit formula for the fundamental solution of the (Kohn) sub-Laplacian on $H(3)$ (which is up to a constant given by $\lVert \mathbf{c} \rVert_{16}^{-2}$). 
In this context we thus obtain an estimate of the fundamental solution of $\mathcal{L}$ by estimating it with $\lVert \mathbf{c} \rVert_{16}^{-2}$, which is proportional to the exact fundamental solution of $\mathcal{L}_0$ on our approximated group $(SE(2))_0$.

\subsubsection{Approximation of the sub-Riemannian distance}
Finally we arrive at the sub-Riemannian distance approximations.
By the Ball-Box theorem (see e.g. \cite{Bellaiche1996}) and equivalence of homogeneous norms, there exists a constant $\gothic{c}$ such that
$$
\gothic{c}^{-1} \lVert \operatorname{Log}(g) \rVert_\zeta \leq d_0(e, g ) \leq \gothic{c} \lVert \operatorname{Log}(g) \rVert_\zeta,
$$
with $\operatorname{Log}(g)$ defined by Eq.~(\ref{eq:coordinatesfirstkind}). The logarithmic norm is locally equivalent to the sub-Riemannian distance, which was proved in full generality in \cite[Thm.~2 \& 4]{nagel-stein-wainger}. Via a scaling of the generators $\tilde{{A}}_2 = \xi^{-1} {A}_2$ and $\tilde{{A}}_3 = \xi^{-1} {A}_3$ we define the $\xi$-isotropic norm
\begin{equation}\label{eq:approxSE2}
\begin{array}{rl}
\lVert \mathbf{c} \rVert_{\xi,\zeta} :&= \sqrt[4]{(|c^1|^2 + |\tilde{c}^2|^2)^2 + \zeta |\tilde{c}^3|^2}\\
&= \sqrt[4]{(|c^1|^2 + \xi^2 |c^2|^2)^2 + \zeta\; \xi^2 |c^3|^2},
\end{array}
\end{equation}
with $\tilde{c}^2 = \xi c^2$ and $\tilde{c}^3 = \xi c^3$, and the $c^i$ given in (\ref{eq:coordinatesfirstkind}). The norm $\lVert \cdot \rVert_{\xi,\zeta}$ closely approximates the sub-Riemannian distance $d_{0}(e,\cdot)$ for $C=1$ (no data-adaptivity) via
\begin{equation}\label{eq:approxDistSE2}
d_{0}(g,h) \approx |\operatorname{Log}(g^{-1}h)|_{\xi,\zeta}, \;\;\;\; |\operatorname{Log}(g)|_{\xi,\zeta}:=\lVert \mathbf{c} \rVert_{\xi,\zeta}
\end{equation}
with $\mathbf{c}$ the coordinates of the first kind obtained via (\ref{eq:coordinatesfirstkind}). In view of the Folland-Kaplan-Kor{\'a}nyi gauge setting $\zeta=16$ in $\lVert \cdot \rVert_{\xi,\zeta}$ would be a sensible choice. We do observe however that $\zeta=44$ gives an even sharper approximation, see Fig.~\ref{fig:DistancesSE2} for a visual comparison to the sub-Riemannian distance $d_0$ and Appendix~\ref{app:optimizationZeta} for a quantitative comparison. The setting $\zeta=44$ is used in all experiments on $SE(2)$.

\section{Sub-Riemannian Distance and its Approximation in $SE(3)$}\label{sec:SE3}
In this section we extend the concepts of the previous section to the group $SE(3)$ of 3D translations and rotations. In the end we again obtain an approximation for the sub-Riemannian distance, which allows us to do perceptual grouping in 3D images as well.

\subsection{The Lie Group $SE(3)$}
\subsubsection{$SE(3)$}\label{subsec:SE3}
The Lie group $SE(3) = \mathbb{R}^3 \rtimes SO(3)$ is the semi-direct product of the group of 3D translations $\mathbb{R}^3$ and the group of 3D rotations $SO(3)$. The group product and inverse for elements $g = (\mathbf{x},\mathbf{R}),g'=(\mathbf{x}',\mathbf{R}') \in SE(3)$ are defined by

\begin{equation}
\begin{array}{rl}
g \cdot g' &= (\mathbf{x}, \mathbf{R}) \cdot (\mathbf{x}',\mathbf{R}') = (\mathbf{x} + \mathbf{R} \mathbf{x}', \mathbf{R} \mathbf{R}'), \\
g^{-1} &= (-\mathbf{R}^{-1}\mathbf{x},\mathbf{R}^{-1}).
\end{array}
\end{equation}
In the 3D case, we define the space of coupled positions and orientations as a Lie group quotient of $SE(3)$:
\begin{equation*}
\mathbb{R}^3 \rtimes S^2 := SE(3)/({\mathbf{0}}\times SO(2)).
\end{equation*}
The group action of $g \in SE(3)$ onto $(\mathbf{y},\mathbf{n}) \in \mathbb{R}^3 \times S^2$ is defined by
\begin{equation*}
g \cdot (\mathbf{y},\mathbf{n}) = (\mathbf{x},\mathbf{R}) \cdot (\mathbf{y},\mathbf{n}) = (\mathbf{x} + \mathbf{R} \mathbf{y}, \mathbf{R} \mathbf{n}).
\end{equation*}
We can identify the element $(\mathbf{x},\mathbf{n}) \in \mathbb{R}^3 \times S^2$ with group elements $(\mathbf{x},\mathbf{R}_\mathbf{n}) \in SE(3)/(\mathbf{0}\times SO(2))$, where $\mathbf{R}_\mathbf{n}$ is any rotation matrix such that $\mathbf{R}_\mathbf{n} \mathbf{e}_z = \mathbf{n}$.

\subsubsection{The Lie Algebra, Exponential Map and Commutators}
Analogously as in the $SE(2)$ case, we associate with the group $SE(3)$ the Lie algebra $\gothic{se}(3)$ using the exponential and logarithmic maps. This is most easily done using an isomorphism with the corresponding matrix group:

\begin{equation*}
(\mathbf{x},\mathbf{R}_{\gamma,\beta,\alpha}) \leftrightarrow \begin{pmatrix}
\mathbf{R}_{\gamma,\beta,\alpha} & \mathbf{x}^T \\
0 & 1
\end{pmatrix}.
\end{equation*}
A basis for the corresponding matrix Lie-algebra is given by
\begin{equation}
\begin{aligned}
\mathbf{X}_1 &= \begin{pmatrix}
0 & 0 & 0 & 1 \\
0 & 0 & 0 & 0 \\
0 & 0 & 0 & 0 \\
0 & 0 & 0 & 0
\end{pmatrix}, \qquad \mathbf{X}_2 = \begin{pmatrix}
0 & 0 & 0 & 0 \\
0 & 0 & 0 & 1 \\
0 & 0 & 0 & 0 \\
0 & 0 & 0 & 0
\end{pmatrix}, \\ \mathbf{X}_3 &= \begin{pmatrix}
0 & 0 & 0 & 0 \\
0 & 0 & 0 & 0 \\
0 & 0 & 0 & 1 \\
0 & 0 & 0 & 0
\end{pmatrix}, \qquad
\mathbf{X}_4 = \begin{pmatrix}
0 & 0 & 0 & 0 \\
0 & 0 & -1 & 0 \\
0 & 1 & 0 & 0 \\
0 & 0 & 0 & 0
\end{pmatrix}, \\ \mathbf{X}_5 &= \begin{pmatrix}
0 & 0 & 1 & 0 \\
0 & 0 & 0 & 0 \\
-1 & 0 & 0 & 0 \\
0 & 0 & 0 & 0
\end{pmatrix}, \qquad \!\!\!\!\mathbf{X}_6 = \begin{pmatrix}
0 & -1 & 0 & 0 \\
1 & 0 & 0 & 0 \\
0 & 0 & 0 & 0 \\
0 & 0 & 0 & 0
\end{pmatrix},
\end{aligned}
\end{equation}
and their equivalents $A_i$ in the tangent space of $SE(3)$ span the Lie algebra $\gothic{se}(3)$. Since it will be clear from the context if we are in the $SE(2)$ or $SE(3)$ case, we use the same notation for the generators of the Lie algebra as previously. Now the left-invariant vector fields are again obtained using the push-forward of the left-multiplication $(L_g)^*$, but they depend on the choice of coordinates. In this paper we mostly rely on $ZYZ$-Euler angles in the parameterization of $SO(3)$, i.e.,
\begin{equation}\label{eq:eulerangles}
\mathbf{R}_{\gamma,\beta,\alpha} = \mathbf{R}_{\mathbf{e}_z,\gamma}  \mathbf{R}_{\mathbf{e}_y,\beta}  \mathbf{R}_{\mathbf{e}_z,\alpha},
\end{equation}
with $\mathbf{R}_{\mathbf{n},\alpha}$ a rotation with angle $\alpha$ around $\mathbf{n}$. Then, the left-invariant vector fields are given by

\begin{equation}
\begin{array}{l}
\mathcal{A}_1|_g = (\cos \alpha \cos \beta \cos \gamma - \sin \alpha \sin \gamma) \partial_x  \vspace{1mm} \\
+ (\sin \alpha \cos \gamma + \cos \alpha \cos \beta \sin \gamma) \partial_y - \cos \alpha \sin \beta \partial_z \vspace{2mm} \\
\mathcal{A}_2|_g = (- \sin \alpha \cos \beta \cos \gamma - \cos \alpha \sin \gamma) \partial_x \vspace{1mm} \\
+ (\cos \alpha \cos \gamma - \sin \alpha \cos \beta \sin \gamma) \partial_y + \sin \alpha \sin \beta \partial_z, \vspace{2mm} \\
\mathcal{A}_3|_g = \sin \beta \cos \gamma \partial_x + \sin \beta \sin \gamma \partial_y + \cos \beta \partial_z \vspace{1mm} \\
\mathcal{A}_4|_g = \cos \alpha \cot \beta \partial_\alpha + \sin \alpha \partial_\beta - \dfrac{\cos \alpha}{\sin \beta}\partial_\gamma, \vspace{1mm} \\
\mathcal{A}_5|_g = -\sin \alpha \cot \beta \partial_\alpha + \cos \alpha \partial_\beta + \dfrac{\sin \alpha}{\sin \beta}\partial_\gamma \vspace{1mm} \\
\mathcal{A}_6|_g = \partial_\alpha,
\end{array}
\end{equation}
for $\beta \neq 0,\pi$.

\begin{remark}
A second coordinate chart is needed to cover the entire $SO(3)$, for which for example $ZYX$-angles can be used, as is done in e.g. \cite{duits2011HARDI}, where also the expressions for the vector fields in this alternative coordinate chart are given. In fact, the basis elements $A_i$ of the Lie algebra correspond to partial derivatives with respect to the $ZYX$-angles, similar to the $SE(2)$-case.
\end{remark}
We can express each element $\gothic{se}(3)$ in terms of the basis with coefficients $\mathbf{c} = (c^1, \dots, c^6)^T$. Furthermore, we define $\mathbf{c}^{(1)}:= (c^1, c^2, c^3)^T$ and $\mathbf{c}^{(2)} := (c^4,c^5,c^6)^T$, the spatial and rotational coefficients, respectively. We can make the exponential map $\operatorname{Exp}_{SE(3)} : \gothic{se}(3) \rightarrow SE(3)$ and logarithmic map $\operatorname{Log}_{SE(3)} : SE(3)\rightarrow \gothic{se}(3)$ explicit using these coefficients. For a $3 \times 3$ matrix $\mathbf{\Omega}$ of the form
\begin{equation}
\mathbf{\Omega} := \left(\begin{array}{ccc}
0 & -c^6 & c^5 \\
c^6 & 0 & -c^4 \\
-c^5 & c^4 & 0
\end{array}
\right),
\end{equation}
we obtain a rotation using the exponential map of matrices, i.e., $\mathbf{R} = \exp(\mathbf{\Omega})$. The relation between the spatial coefficients $\mathbf{c}^{(1)}$ and $(\mathbf{x},\mathbf{R})$ is given by

\begin{equation}\label{eq:c1}
\mathbf{c}^{(1)}  = \left(I - \frac12 \mathbf{\Omega} + q^{-2} \left(1 - \frac{q}{2} \cot \left(\frac{q}{2} \right) \right)(\mathbf{\Omega})^2 \right) \mathbf{x},
\end{equation}
where $q = ||\mathbf{c}^{(2)}||$ and $\mathbf{\Omega}$ such that $\mathbf{R} = \exp(\mathbf{\Omega})$. Now

\begin{equation}
\begin{aligned}\label{eq:ExpLogSE3}
&\operatorname{Log}_{SE(3)}(g) = \sum_{i=1}^6 c_i(g) A_i, \qquad \text{and}\\
&\operatorname{Exp}_{SE(3)}\left(\sum_{i=1}^6 c_i(g) A_i \right) = g,
\end{aligned}
\end{equation}
using the relations above.

\subsection{Sub-Riemannian Geometry in SE(3)}\label{subsec:SRGeometrySE3}
In the SE(3) case, a horizontal path is a curve $\gamma:\mathbb{R}\rightarrow SE(3)$ with tangent vectors $\dot{\gamma}(t) \in \left.\Delta\right|_{\gamma(t)}:= \operatorname{span} \{ \left. \mathcal{A}_3\right|_{\gamma(t)}, \left. \mathcal{A}_4 \right|_{\gamma(t)}, \left. \mathcal{A}_5 \right|_{\gamma(t)}\}$, where $\Delta$ is now the sub-bundle of full tangent bundle spanned by $\{\mathcal{A}_i\}_{i=1}^6$. In this case we have one spatial control $u^3$ and two `angular' controls $u^4$ and $u^5$, so that the sub-Riemannian metric tensor becomes:

\begin{equation}\label{eq:SRMetricTensorSE3}
\begin{array}{rl}
  \left.\mathcal{G}^{\xi,C}\right|_{\gamma(t)} (\dot{\gamma}(t),\dot{\gamma}(t)):= &C(\gamma(t))^2 \left(\xi |u^3(t)|^2 +\right.\\ &\left.\quad |u^4(t)|^2 + |u^5(t)|^2\right),
  \end{array}
\end{equation}

The sub-Riemannian distance between two elements $g_1,g_2 \in SE(3)$ is still defined as in (\ref{eq:SRDistancesSE2}), but now the infimum is taken over Lipschitz continuous curves $\gamma \in \operatorname{Lip}([0,T],SE(3))$ with $\gamma(0) = g_1$, $\gamma(1) = g_2$, and $\dot{\gamma}(t) = u^3(t) \left.\mathcal{A}_3\right|_{\gamma(t)} + u^4(t) \left.\mathcal{A}_4\right|_{\gamma(t)} + u^5(t) \left.\mathcal{A}_5\right|_{\gamma(t)}$.

\subsection{A Nilpotent Approximation $(SE(3))_0$ of $SE(3)$ and the Approximated sub-Riemannian Distance}\label{subsec:approxSE3}

It is important to realize that the logarithmic map is only well-defined on the group $SE(3)$ and not on the quotient $\mathbb{R}^3 \rtimes S^2$, i.e., different choices for $\alpha$ in the rotational part result in different values for the coefficients $c^i$. Here, we choose the approach of \cite{portegies2015ssvm} and set $\alpha=-\gamma$ such that expected symmetries are preserved. With that choice, the logarithm (\ref{eq:ExpLogSE3}) gives for each $(\mathbf{x},\mathbf{n}) \in \mathbb{R} \rtimes S^2$ a unique vector $\mathbf{c}$, on which we can put a norm:

\begin{equation}
\begin{aligned}\label{eq:FKKGaugeSE3}
&|\operatorname{Log}_{SE(3)}(g)|_{\xi,\zeta} := ||\mathbf{c}||_{\xi,\zeta} :=\\
&\sqrt[4]{ (\xi^2 |c^3|^2 + |c^4|^2 + |c^5|^2)^2 +   \zeta \, (\xi^2(|c^1|^2 + |c^2|^2) + |c^6|^2) },
\end{aligned}
\end{equation}
where $\mathbf{c} = \mathbf{c}(g)$ according to (\ref{eq:ExpLogSE3}).

Also here, the Folland-Kaplan-Kor{\'a}nyi-type norm can be used to approximate the fundamental solutions of the sub-Laplacian on SE(3). The norm $||\mathbf{c}||_{\xi,\zeta}$ with $\zeta=1$ was for example used in \cite{duits2011HARDI} approximations of the heat kernel and the fundamental solution on $SE(3)$, of which only recently exact solutions were found in \cite{portegies2016arxiv}. In the context of this paper, we can approximate the exact solutions of the sub-Laplacian on $SE(3)$ by $\lVert \textbf{c} \rVert_{1,16}^{2-Q}$, with homogeneous dimensions $Q=9$, as the exact solution of the sub-Laplacian on the approximation group $(SE(3))_0$. The group $(SE(3))_0$ that locally approximates $SE(3)$, is again obtained via a nilpotent step 2 approximation of the BCH formula, and is defined by the group product
\begin{multline} \label{eq:GroupLawNilpotentSE3}
 (x^1,x^2,x^3,x^4,x^5,x^6) \cdot (y^1,y^2,y^3,y^4,y^5,y^6) \\=
 \left(\begin{array}{c}
 x^1 + y^1 + \frac{1}{2}(x^5y^3 - x^3y^5)\\
 x^2 + y^2 + \frac{1}{2}(x^3y^4 - x^4y^3) \\
 x^3 + y^3 \\
 x^4 + y^4 \\
 x^5 + y^5 \\
 x^6 + y^6 + \frac{1}{2}(x^4y^5 - x^5y^4) \\
 \end{array}\right)^T,
\end{multline}
with $x^i,y^i$ coordinates of the first kind given by the logarithmic map (\ref{eq:c1}). This new group is a free-nilpotent group of rank 3 and step 2.

\begin{figure*}
  \centering
  \includegraphics[width=\textwidth]{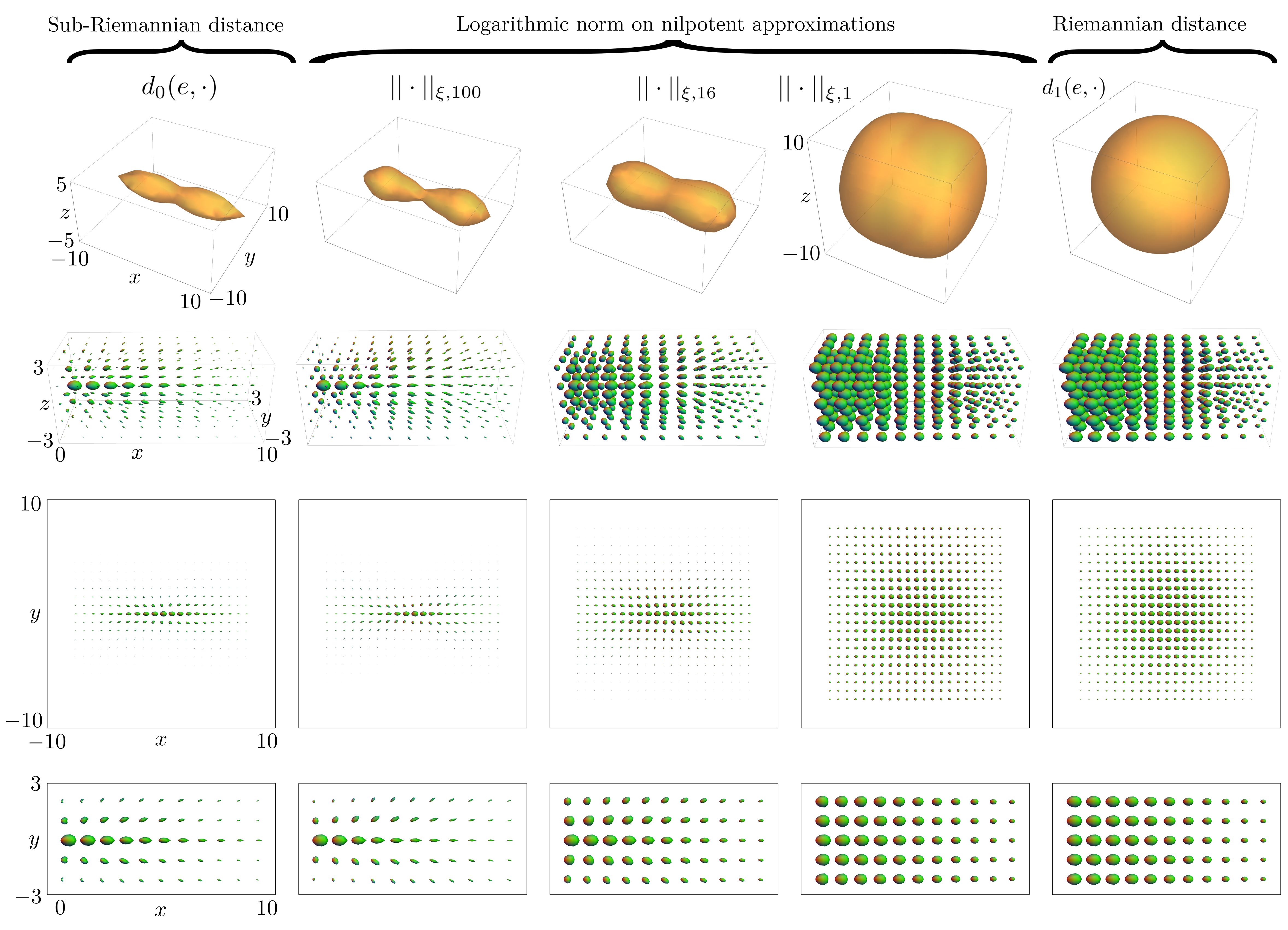}
  \caption{Distances on $SE(3)$ for $\xi=.1$, $C=1$, with the origin placed at $e=(\mathbf{0},\mathbf{e}_x)$. Top row: Level sets of the spatial projections (minimum intensity projections over $S^2$) of the distance volumes on $SE(3)$. Rows two to four: Glyph visualizations in which each distance volume $d$ is visualized with a "Gaussian" density $U(g)=e^{-d(e,g)^2}$. For an interpretation of the glyphs see Remark~\ref{rem:glyph}. Row two: Glyph visualizations of sub-volume. Row three: Glyph visualization of slice at $z=0$. Row four: Zoomed in glyph visualization of the slice a $z=0$. From left to right: The sub-Riemannian distance $d_{0}(e,\cdot)$ on $SE(3)$, see Eq.~(\ref{eq:SRDistancesSE2}) and (\ref{eq:SRMetricTensorSE3}); Homogenous norms $\lVert \cdot \rVert_{\xi,\zeta}$, see Eq.~(\ref{eq:FKKGaugeSE3}), of the nilpotent approximation $(SE(3))_0$ for respectively $\zeta=100$, $\zeta=16$ (Folland-Kaplan-Kor{\'a}nyi gauge) and $\zeta=1$; The ($\xi$-isotropic) Riemannian distance $d_1(e,\cdot)$ on $SE(3)$, see Table~\ref{tab:metrics} for an overview of the different distances. }\label{fig:DistancesSE3}
\end{figure*}

We approximate the sub-Riemannian distance $d_{0}$ on $SE(3)$ via the norm (\ref{eq:FKKGaugeSE3}). I.e.,
\begin{equation}\label{eq:approxDistanceSE3}
d_{0}(g,h) \approx |\operatorname{Log}_{SE(3)}(g^{-1}h)|_{\xi,\zeta},
\end{equation}
 and as such again obtain an approximation of the distance in the sense of Rothschild and Stein \cite{rothschild-stein}.
Based on the quantitative comparison to the sub-Riemannian distances $d_0$ in Appendix \ref{app:optimizationZeta} and the visualizations in Fig.~\ref{fig:DistancesSE3} of the level sets we conclude that the approximated sub-Riemannian distance of (\ref{eq:approxDistanceSE3}) quite accurately approximates the true sub-Riemannian distance on $SE(3)$. In our analysis we found that the logarithmic norm with $\zeta=100$ gave the best approximation, and as such we used this norm in the perceptual grouping experiments of Sec.~\ref{subsec:PGSE3}.

\begin{remark}\label{rem:glyph}
The glyph at each grid point $\mathbf{y}$ in Fig.~\ref{fig:DistancesSE3} is given by the surface $\{\mathbf{y} + \nu U(\mathbf{y}, \mathbf{n})\mathbf{n} | \mathbf{n}\in S^2\}$, for a specific choice $\nu > 0$, and with density $U:\mathbb{R}^3\times S^2\rightarrow\mathbb{R}^+$. The color of each orientation $\mathbf{n} = (n^1,n^2,n^3) \in S^2$ on the glyph is defined by the RGB color $(n^1,n^2,n^3)$.
\end{remark}

\section{Perceptual Grouping, Fast Marching and Key Point Tracking}\label{sec:algorithms}
In this section the algorithms used in this paper are explained. Our main application of interest is that of grouping/clustering of points on blood vessels via the perceptual grouping algorithm, which is explained in Subsec.~\ref{subsec:PG}. The perceptual grouping algorithm takes as input a set of key points that
are obtained via the \emph{minimal path tracking with key points} algorithm \cite{benmansour2009fast}, explained in Subsec.~\ref{subsec:KP}, which is an adaptation of the fast marching algorithm, explained in Subsec.~\ref{subsec:FM}. Finally since different metrics are used throughout the experiments (both for generating key points and for perceptual grouping) we end this section with an overview of the used metrics in this paper in Subsec.~\ref{subsec:metrics}.

\subsection{The Perceptual Grouping Algorithm}\label{subsec:PG}
The perceptual grouping algorithm presented in this paper is a modification of the original algorithm proposed by Cohen \cite{Cohen2001multiple}, and which was later adapted for perceptual grouping based on anisotropic distances \cite{Bougleux2008anisotropic}. In recent work \cite{Chen2017}, the perceptual grouping algorithm was extended for the grouping of $n$ closed contours for an a-priori specified $n$. Like in \cite{Bougleux2008anisotropic} and \cite{Chen2017}, we use the main algorithm of \cite{Cohen2001multiple} as a backbone, but we change the metric used for perceptual grouping and we impose an additional constraint to avoid closed loops (which are physically not realistic in the vessel networks of interest). Our adapted perceptual grouping algorithm is given in pseudo code in Algorithm~\ref{alg:PerceptualGrouping}.

The goal of the perceptual grouping algorithm is to construct a graph out of a set $\mathcal{S}$ of points of interest in which the edges $\mathcal{D}_\mathcal{S}$ are true connections (represented by geodesics) between points. Following the terminology of \cite{deschamps2001fast,benmansour2009fast,chen2016vessel} we will refer to the points of interest as \emph{key points}. Each key point is only linked to at most 2 other key points (i.e., node degree $\delta_i$ is 2 at most). The final graph thus only contains sets of non-bifurcating vessel segments. The graph is build up by inserting one-by-one the edges which have the shortest geodesic distance (if the node degree allows). As such, only the strongest connections (shortest geodesics) appear in the final graph network. Since the original algorithm in \cite{Cohen2001multiple} (and also \cite{Bougleux2008anisotropic}) does not include a mechanism to avoid closed loops we include an additional check in the main algorithm to prevent this. Finally, in order to avoid connecting key points which are too far apart from each other we only consider edges of which the spatial arc length of the connecting geodesic does not exceed a certain a-priori threshold $s_{max}$.

In summary, our changes relative to the works \cite{Cohen2001multiple,Bougleux2008anisotropic,Chen2017}, is that we
\begin{itemize}
  \item keep the choice for distance $d(x_i,x_j)$ open. In our experiments the distances $d$ will be mainly based on sub-Riemannian geometry in $SE(n)$.
  \item explicitly avoid making long distance connections by filtering out such possible connections in an initialization step.
  \item avoid closed loops by not making connections between nodes that belong to the same sub-graph.
  \item group \emph{crossing lines} without pre-specifying the number of groups.
\end{itemize}
In particular, it is the use of a sub-Riemannian metric on $SE(n)$ that allows for the grouping of crossing lines. A first (successful) feasibility study on the possibility of perceptual grouping of crossing lines was performed by Chen et al. \cite{Chen2017} using a (sub-)Finsler metric (based on the Euler elastica model) on position-orientation space. There it was successfully demonstrated on phantom images that their algorithm is able to deal with crossing closed contours, however, it required specification of the number of contours (which is not always a-priori known). Furthermore, their metric relies on a notion of directionality (instead of just orientations) which is useful in grouping closed contours, but may be disadvantages for grouping non-closed contours. Here, we focus on the grouping of non-closed crossing contours without specifying the number of contours. Furthermore, we quantify the performance of perceptual grouping of crossing lines on a large set of both retinal images in 2D, and phantom images in 3D.



\begin{algorithm}[h]\label{alg:PerceptualGrouping}
\hrulefill\\
\footnotesize	
\SetKwData{Left}{left}\SetKwData{This}{this}\SetKwData{Up}{up}
\SetKwFunction{Union}{Union}\SetKwFunction{FindCompress}{FindCompress}
\SetKwInOut{Input}{input}\SetKwInOut{Output}{output}\SetKwInOut{Variables}{variables}
 \Input{
    $\mathcal{S}$: a set of key points\; \\
 	\hspace{0.3em}$d(g_i,g_j)$: distances between $g_i, g_j \in \mathcal{S}$\; \\
    \hspace{0.3em}$s_{max}$: max spatial length of geodesics\;
    }
 \Variables{
    $\tilde{\mathcal{D}}_\mathcal{S}$: set of possible edges\; \\
    \hspace{0.3em}$\delta_i$: node degree of $x_i$\;
    }
 \Output{
    $\mathcal{D}_\mathcal{S}$: final set of edges\;
	}
\vspace{1em}
\emph{Initialization}: \newline
Compute the distances $d(g_i,g_j)$ (and corresponding geodesics) between all key points $g_i,g_j\in\mathcal{S}$.\vspace{0.5em}\newline
Initialize $\tilde{\mathcal{D}}_\mathcal{S}$ with the set of all edges between each $g_i,g_j \in \mathcal{S}$ whose connecting geodesic has spatial arc length smaller then $s_{max}$, and set $\mathcal{D}_\mathcal{S}=\O$.\\
\vspace{1em}
\emph{Main algorithm}: \\
\While{$\tilde{\mathcal{D}}_\mathcal{S} \neq \O$ }{
\begin{enumerate}
  \item  \emph{Select edge and remove it from $\tilde{\mathcal{D}}_\mathcal{S}$}: \\
  \hspace{1em} $(g_i,g_j) \leftarrow \underset{(g,h)\in\tilde{\mathcal{D}}_\mathcal{S}}{\operatorname{argmin}} d(g,h)$\;
  \hspace{1em} $\tilde{\mathcal{D}}_\mathcal{S} = \tilde{\mathcal{D}}_\mathcal{S} - (x_i,x_j)$\;
  \item \emph{Check topology and update network}: \\
  \hspace{1em} {\bf if} $\delta_i < 2$ and $\delta_j < 2$ {\bf and} $g_i,g_i$ are not\\
  \hspace{1em} already in the same sub-graph in $\mathcal{D}_\mathcal{S}$ \\
  \hspace{1em} {\bf then} $\mathcal{D}_\mathcal{S} = \mathcal{D}_\mathcal{S} + (g_i,g_j)$\;$\delta_i=\delta_i+1$\;
\end{enumerate}
  }
  \hrulefill
 \caption{Perceptual grouping.}
\end{algorithm}

\subsection{Fast Marching}\label{subsec:FM}
Most of the distances (except for the fast analytic approximations) and the geodesics used in this paper are computed via the fast marching algorithm, which is an efficient numerical solver of the eikonal equation and which can be used to obtain \emph{globally optimal} geodesics \cite{cohen1997global}. Let $g_0$ be an arbitrary source point in a domain $\mathbb{M}$ of interest,
let $\mathcal{G}|_{g}: T_g(\mathbb{M}) \times T_g(\mathbb{M}) \rightarrow \mathbb{R}^+$ be a metric tensor defined on the tangent space $T_g(\mathbb{M})$ at $g\in\mathbb{M}$, and let
\begin{equation}\label{eq:distancemap}
U(g)\hspace{-0.25em}:=\hspace{-0.25em}d(g_0,g) \hspace{-0.25em}= \hspace{-0.45em}\underset{\gamma \in \mathcal{S}(g_0,g)}{\operatorname{inf}} \int_0^1 \sqrt{ \left.\mathcal{G}\right|_{\gamma(t)}(\dot{\gamma}(t),\dot{\gamma}(t)) }{\rm d}t
\end{equation}
its associated distance map, where the infimum is taken over the set $\mathcal{S}(g_0,g)$ of Lipschitz continuous curves with $\gamma(0)=g_0$, $\gamma(1) = g$, and with $\dot{\gamma}(t)\in T_{\gamma(t)}(\mathbb{M})$. Then the distance map $U$ is the unique viscosity solution of the eikonal equation
\begin{equation}\label{eq:eikonal}
  \left\{
  \begin{array}{l}
  \sqrt{\mathcal{G}\left( \nabla_\mathcal{G} U(g),\nabla_\mathcal{G} U(g) \right)} = 1,\\
  U(g_0) = 0,
  \end{array}
  \right.
  \leftrightarrow
  \left\{
  \begin{array}{l}
  \lVert \nabla_\mathcal{G} U(g) \rVert_\mathcal{G} = 1,\\
  U(g_0) = 0,
  \end{array}
  \right.
\end{equation}
with $\nabla_\mathcal{G} := \mathcal{G}^{-1} {\rm d}U$ the intrinsic gradient with inverse metric $\mathcal{G}^{-1}$  and ${\rm d}U$ the differential of $U$, and $\lVert \cdot \rVert_\mathcal{G}$ the norm with respect to the metric tensor. In the standard (data-adaptive) Euclidean case with $\mathbb{M}=\mathbb{R}^2$, $g_0=\mathbf{0}$, $g=\mathbf{x}$, $\dot{\gamma}(t) = u^1(t) \partial_x + u^2(t) \partial_y \in T_{\gamma(t)}(\mathbb{R}^2)$, and with $\mathcal{G}|_{\gamma(t)}(\dot{\gamma}(t),\dot{\gamma}(t)) = C(\gamma(t))^2 (|u^1(t)|^2+|u^2(t)|^2)$ the eikonal equation is given by $\lVert \nabla U (\mathbf{x})\rVert = C(\mathbf{x})$.

The fast marching algorithm efficiently solves the eikonal equation in a one pass algorithm. It computes the values of $U$ in increasing order (starting with $U(g_0)=0$) based on the Bellman principle of optimality, in a manner very similar to the Dijkstra algorithm for shortest paths on graphs \cite{dijkstra1959note}. The minimal geodesic connecting $g_0$ with $g$ is then obtained via a gradient descent on $U$ from $g$ back to the origin $g_0$, i.e., solving the ODE
\begin{equation*}
\left\{
\begin{array}{l}
\dot{\gamma}(t) \propto -\mathcal{G}^{-1} {\rm d}U(\gamma(t)),\\
\gamma(0) = g_0.
\end{array}
\right.
\end{equation*}
For details on the fast marching algorithm on isotropic manifolds we refer to \cite{Tsitsiklis1995,Sethian1999}, to \cite{Mirebeau2014,Jbabdi2008} for anisotropic fast marching, and to \cite{Sanguinetti2015CIARP} and \cite{duits2016optimal} for fast marching in sub-Riemannian manifolds in $SE(2)$ and $SE(3)$ respectively.


\subsection{Generating Key Points}\label{subsec:KP}
The key point method is based on keeping track of a spatial arc-length map $U_l$ (in which the spatial lengths of the minimizing geodesics $\gamma$ defining $U$ are stored), and stop as soon as a certain distance threshold is passed \cite{deschamps2001fast}. The rationale behind this algorithm is that among all points with equal geodesic distance values $U$, the points reached by geodesics $\gamma$ that best follow the data (paths along which $C$ is low) have maximum spatial distance $l(\gamma)$. Such a point maximizing spatial distance in a given level set in $U$ is called a key point. The fast marching algorithm is highly suited for keeping track of a spatial arc-length map $U_l$, in addition to $U$, due to the local updating approach (wavefront propagation). Moreover, the algorithm can stop early if one is only interested in finding the first key point with length larger than $l_{max}$ \cite{deschamps2001fast}.

In summary, a key point is detected as follows. The spatial arc-length map is defined as
\begin{equation}\label{eq:spatialdistancemap}
U_l(g):= l(\gamma_{g_0,g}),
\end{equation}
with $\gamma_{g_0,g}= \underset{\gamma \in \mathcal{S}(g_0,g)}{\operatorname{argmin}} \int_0^1 \sqrt{ \left.\mathcal{G}\right|_{\gamma(t)}(\dot{\gamma}(t),\dot{\gamma}(t)) }{\rm d}t$ the minimizing geodesic in (\ref{eq:distancemap}), and with
\begin{equation}
l(\gamma) = \int_0^1 \lVert \dot{\mathbf{x}}(t) \rVert {\rm d}t
\end{equation}
the spatial arc-length of $\gamma$, with $\dot{\mathbf{x}}(t) = \mathbb{P}_{\mathbb{R}^n} \dot{\gamma}(t) \in \mathbb{R}^n$ the spatial components of the tangents $\dot{\gamma}(t)$\footnote{In the lifted problem $SE(2)$ the spatial components are for example given by $\dot{\mathbf{x}}(t) = u^2(t) \mathcal{A}_2|_{\gamma(t)} + u^3(t) \mathcal{A}_3|_{\gamma(t)}$, and in the $SE(3)$ case $\dot{\mathbf{x}}(t) = u^1(t) \mathcal{A}_1|_{\gamma(t)} + u^2(t) \mathcal{A}_2|_{\gamma(t)} + u^3(t) \mathcal{A}_3|_{\gamma(t)}$.}. The fast marching algorithm stops as soon as there is a $g$ for which $U_l(g)\ge l_{max}$, and $g$ will be called a key point.

With the above criteria one can iteratively detect new key points based on geodesic distances to previously found key points, a method known as \emph{minimal path tracking with key point detection} \cite{benmansour2009fast}. One can make several choice on when to stop the key point tracking algorithm \cite{benmansour2009fast,kaul2012detecting,chen2016vessel}. In this work, we rely on the approach by Chen et al. \cite{chen2016vessel}, where we only add key points on locations which lie in a masked region (we use a binary vessel centerline mask $m:\mathbb{M}\rightarrow\{0,1\}$), i.e., we only add a key point when both $U_l(g)\ge l_{max}$ and $m(g)=1$. The algorithm is stopped as soon $U_l(g)\ge 3 \; l_{max}$.


\subsection{Overview of Distances Used in This Paper}\label{subsec:metrics}
Table~\ref{tab:metrics} gives an overview of the different distances discussed in this paper and used in the experiments. The isotropic Euclidean metrics are used the generate key points in $\mathbb{R}^2$ and $\mathbb{R}^3$ using the algorithm of Subsec.~\ref{subsec:KP}. The isotropic Euclidean distances are also used in comparison to the other distances on $SE(n)$ in the perceptual grouping experiments. The sub-Riemannian distances on $SE(2)$ and $SE(3)$ are explained respectively in Subsec.~\ref{subsec:SRGeometrySE2} and Subsec.~\ref{subsec:SRGeometrySE3}. In the Riemannian distances the full tangent bundle on $SE(n)$ is considered. This means that now also non-horizontal curves in $SE(n)$ are considered, i.e., points on the curves $\gamma$ are allowed to move sideways by the non-horizontal controls $u^3(t)$ in the $SE(n)$ case, and $u^1(t), u^2(t)$ in the $SE(3)$ case. Recall that in this case the blue and red oriented particles in Fig.~\ref{fig:AlignmentSE2} do have the same distance to the source (black arrow). Finally, the sub-Riemannian distance approximations, denoted with $|\operatorname{Log}_{SE(n)}(g^{-1}h)|_{\xi,\zeta}$
are discussed and defined in respectively Subsec.~\ref{sec:ApproximationSE2} and Eq.~(\ref{eq:approxDistSE2}) for $SE(2)$ and  Subsec.~\ref{subsec:approxSE3} and Eq.~(\ref{eq:FKKGaugeSE3}) for $SE(3)$.

\begin{table*}
\caption{Overview of the metrics used in this paper.}
\centering
\begin{tabular}{lllll}
\toprule
Distance & Manifold       & Tangent b.  & Tangent vectors $\dot{\gamma}$ & Metric tensor  $\mathcal{G}$\\
notation & $\mathbb{M}$   & $T(\mathbb{M})$ &   &  \\
\midrule
\multicolumn{5}{c}{Isotropic Euclidean$^*$}\\
\midrule
$\lVert g-h \rVert$  & $\mathbb{R}^2$ & $T(\mathbb{R}^2)$ & $\dot{\gamma}(t) = u^1(t) \partial_x + u^2(t) \partial_y$ & $\mathcal{G}|_{\gamma(t)} = C(\gamma(t))^2(|u^1(t)|^2 + |u^2(t)|^2)$ \\
\\
$\lVert g-h \rVert$  & $\mathbb{R}^3$ & $T(\mathbb{R}^3)$
& $\dot{\gamma}(t) = \begin{array}{l}u^1(t) \partial_x + u^2(t) \partial_y\\+ u^3(t) \partial_z\end{array}$
& $\mathcal{G}|_{\gamma(t)} = C(\gamma(t))^2\left(\begin{array}{l}|u^1(t)|^2 + |u^2(t)|^2\\ + |u^3(t)|^2\end{array}\right)$ \\
\midrule
\multicolumn{5}{c}{(Full) Riemannian $SE(n)$}\\
\midrule
$d_{1}(g,h)$  & $SE(2)$ & $T(SE(2))$
& $\dot{\gamma}(t) = \sum\limits_{i=1}^3 u^i(t) \mathcal{A}_i|_{\gamma(t)}$
& $\mathcal{G}|_{\gamma(t)} = C(\gamma(t))^2\left(\sum\limits_{i=1}^3 \xi_i^2 \, |u^i(t)|^2\right)$, \\
& & & & with $\xi_2=\xi_3=\xi$ and $\xi_1=1$\\
\\
$d_{1}(g,h)$  & $SE(3)$ & $T(SE(3))$
& $\dot{\gamma}(t) = \sum\limits_{i=1}^5 u^i(t) \mathcal{A}_i|_{\gamma(t)}$
& $\mathcal{G}|_{\gamma(t)} = C(\gamma(t))^2\left(\sum\limits_{i=1}^5 \xi_i^2 \, |u^i(t)|^2\right)$, \\
& & & & with $\xi_1=\xi_2=\xi_3=\xi$ and $\xi_4=\xi_5=\xi_6=1$\\
\midrule
\multicolumn{5}{c}{Sub-Riemannian $SE(n)$}\\
\midrule
$d_{0}(g,h)$  & $SE(2)$ & $\Delta$
& $\dot{\gamma}(t) = u^1(t) \mathcal{A}_1|_{\gamma(t)} + u^2(t) \mathcal{A}_2|_{\gamma(t)}$
& $\mathcal{G}|_{\gamma(t)} = C(\gamma(t))^2(|u^1(t)|^2 + \xi^2 |u^2(t)|^2)$ \\
\\
$d_{0}(g,h)$  & $SE(3)$ & $\Delta$
& $\dot{\gamma}(t) = \begin{array}{l} u^3(t) \mathcal{A}_3|_{\gamma(t)} + u^4(t) \mathcal{A}_4|_{\gamma(t)}\\+ u^5(t) \mathcal{A}_5|_{\gamma(t)}\end{array}$
& $\mathcal{G}|_{\gamma(t)} = C(\gamma(t))^2\left(\begin{array}{l}\xi^2 |u^3(t)|^2 + |u^4(t)|^2\\ + |u^5(t)|^2\end{array}\right)$ \\
\midrule
\multicolumn{5}{c}{Sub-Riemannian approximation}\\
\midrule
\multicolumn{2}{l}{$|\operatorname{Log}_{SE(2)}(g^{-1}h)|_{\xi,\zeta}$} & \multicolumn{3}{c}{\emph{Approximation of the sub-Riemannian distance on $SE(2)$, cf. Eq.~(\ref{eq:approxDistSE2})}}\\
\multicolumn{2}{l}{$|\operatorname{Log}_{SE(3)}(g^{-1}h)|_{\xi,\zeta}$} & \multicolumn{3}{c}{\emph{Approximation of the sub-Riemannian distance on $SE(3)$, cf. Eq.~(\ref{eq:FKKGaugeSE3})}}\\
\bottomrule
\multicolumn{5}{p{16cm}}{\emph{\small $^*$ The isotropic Euclidean distances are used in key point generation and perceptual grouping. The other distances are only used in the perceptual grouping algorithm.}}
\end{tabular}%
\label{tab:metrics}
\end{table*}

\subsubsection{The Cost $C$}\label{subsubsec:thecost}
The cost functions $C$ are constructed from functions $U_f:\mathbb{R}^n \times S^{n-1}\rightarrow \mathbb{R}$ on the orientation-lifted space. These functions $U_f$ are obtained via an orientation score transform \cite{Duits2007,Janssen2015ssvm} of image $f:\mathbb{R}^n\rightarrow\mathbb{R}$ by correlating the image with a set of anisotropic wavelets $\psi:\mathbb{R}^n \rightarrow \mathbb{R}$:
\begin{equation}
U_f(g) = ( \mathcal{U}_g \psi , f )_{\mathbb{L}_2(\mathbb{R}^n)},
\end{equation}
with $( f , g )_{\mathbb{L}_2(\mathbb{R}^n)}=\int_{\mathbb{R}^n} \overline{f(\mathbf{x})}g(\mathbf{x}) {\rm d}\mathbf{x}$ the standard inner product on $\mathbb{L}_2(\mathbb{R}^n)$, with the overline denoting complex conjugation, and where $\mathcal{U}_g$ denotes the left-regular representation of the Lie group on images $f$. For the group $SE(2)$ acting on images $f \in \mathbb{L}_2(\mathbb{R}^2)$ it is defined as
$$
(\mathcal{U}_g f)(\mathbf{y}):=f(\mathbf{R}_\theta^{-1}( \mathbf{y} - \mathbf{x}) )
$$
with $g = (\mathbf{x},\theta) \in SE(2)$  (recall the group definitions in Subsec.~\ref{subsec:SE2}). For the group $SE(3)$ acting on images $f \in \mathbb{L}_2(\mathbb{R}^3)$ it is defined as
$$
(\mathcal{U}_g f)(\mathbf{y}):=f(\mathbf{R}_{\mathbf{n}}^{-1}( \mathbf{y} - \mathbf{x}) )
$$
with $g = (\mathbf{x}, \mathbf{R}_{\mathbf{n}}) \in SE(3)$ (recall the group definition in Subsec.~\ref{subsec:SE3}).

The wavelets used in the orientation score transform \cite{Duits2007,Janssen2015ssvm} are designed in such a way that all rotated version together cover the full Fourier spectrum. With this design no data is lost in the transformation and a stable invertible transform (from orientation score) back to image exists. For details on this wavelet-type transform for lifting 2D images to functions on $SE(2)$ we refer to \cite{Duits2007}, and for lifting 3D images to 3D orientation scores we refer to \cite{Janssen2015ssvm}. In all experiments we define the cost in the following form
\begin{equation}\label{eq:cost}
C(g) = \frac{1}{1 + \lambda \mathcal{V}(g)^p},
\end{equation}
with $\mathcal{V}$ a vessel (or center line) enhancement obtained by processing of the orientation score $U_f$, and which is normalized between 0 and 1. Parameters $\lambda$ and $p$ then control respectively the influence of the cost (data-adaptivity) and $p$ the contrast.

Good choices for $\mathcal{V}$ for tracking of vessels in 2D position orientation space may be via the vessel enhancements of \cite{Zhang2016} or \cite{Hannink2014}, similar to the $SE(2)$ tracking experiments in \cite{Bekkers2015SIIMS}. For tracking in 3D orientation scores $\mathcal{V}$ may be obtained via the crossing preserving vessel enhancements of \cite{DuitsJanssen2016}. In related tracking problems in lifted spaces the lifts are obtained via tubularity measures \cite{Chen2017,Law2008,LiYezzi}, or by correlating the image with a set of rotated templates \cite{Pechaud}.

\subsubsection{Projective Line Bundle}
Finally, we remark that when dealing with geodesic distances in $SE(n)$ we have to take into account that these are defined for positions and orientations on the \emph{full} sphere $S^{n-1}$. The distances discussed in this paper thus make a distinction between forward and backward arrival directions, i.e., $d(e,(\mathbf{x},\theta)) \neq d(e,(\mathbf{x},\theta + \pi))$.

In practice, and in particular in our perceptual grouping problem, we often do not know the direction of the vessel, but we only have orientations. As such, we would actually want to compute distances on the projective line bundle $\mathbb{R}^n \times P^{n-1}$, with $P^{n-1} := S^{n-1}/\sim$ with identification of antipodal points $\mathbf{n}_1 \sim \mathbf{n}_2 \leftrightarrow \mathbf{n}_1 = \pm \mathbf{n}_2$. We define the distances $\tilde{d}$ on the projective line bundle by distances $d$ on $SE(n)$ via
\begin{equation}
\tilde{d}(g,(\mathbf{x},\mathbf{n})) = \operatorname{min}\left\{ d( g , (\mathbf{x},\mathbf{n})), d( g , (\mathbf{x},-\mathbf{n}) ) \right\},
\end{equation}
with $n \in S^{n-1}$, and $g,(\mathbf{x},\pm\mathbf{n})\in SE(n)$. Note that in the $SE(2)$ case we have with $\mathbf{n}(\theta)= (\cos \theta, \sin \theta) \leftrightarrow \theta$ and $-\mathbf{n}(\theta)=\mathbf{n}(\theta+\pi)$. For a more detailed analysis on data-adaptive sub-Riemannian geodesics on the 2D projective line bundle we refer \cite{bekkers2017vessel}.

\section{Experiments}\label{sec:experiments}
In the experiments we aim to quantify the performance of perceptual grouping with different distances. For a fair comparison we therefore generate automatically the most reasonable key points by using a vessel center line mask $m:\mathbb{R}^n \rightarrow [0,1]$ (see Subsec.~\ref{subsec:KP}) based on the ground truth data. Moreover, this guarantees that the key points are always located on the ground truth center lines, which allows us to quantify performance using the ground truth data. In both the 2D and 3D case the key points are then generated using the isotropic Euclidean metric tensor, and with $\mathcal{V}(\mathbf{x})=m(\mathbf{x})$ (see Subsec.~\ref{subsubsec:thecost}). In all experiments we set $p=1$, $\lambda=100$ to compute the cost (cf.~Eq.~(\ref{eq:cost})).

In the perceptual grouping experiments the cost functions are constructed from orientation score transforms $U_f$ of the mask $m$ on $\mathbb{R}^n$. The costs on $SE(n)$ are then constructed via the modulus of the score:
\begin{equation}\label{eq:vesselnessSEn}
\mathcal{V}(g)=\mathcal{V}_{SE(n)}(g):=|U_f(g)|.
\end{equation}
For equal comparison the costs on $\mathbb{R}^n$ are then constructed via $V(\mathbf{x}) = \underset{\mathbf{n}\in S^{n-1}}{\operatorname{max}}\mathcal{V}_{SE(n)}(\mathbf{x},\mathbf{n})$, i.e., via a maximum intensity projection over orientations $\mathbf{n}$.

\subsection{Perceptual Grouping in $SE(2)$}\label{subsec:PGSE2}
\subsubsection{Experimental Setup}
The data for the 2D retinal vessel grouping experiments consists of 52 retinal image patches in which the vessels have complicated topologies (each patch contains at least 1 crossing, and at least 1 bifurcation). For each retina patch the center lines were semi-automatically traced, after which the connectivity (bifurcation relations) between the vessel segments were manually determined. The set of images contained in total 313 separate vessel segments. A connection between two nodes was determined to be a true positive if both nodes lie on the same vessel tree.

The minimum distance between key points in the retina experiments (with patch sizes of $\approx 400\times400$ pixels) was set to $l_{max}=30$ pixels. The maximum geodesic arc length distance in the perceptual grouping algorithm was set to $s_{max}=80$ pixels. The orientations $\theta$ at each key point $\mathbf{x}$ was estimated by the orientation that gave maximum response in the orientation score, i.e., $\theta=\underset{\theta \in S^1}{\operatorname{argmax}} \; \mathcal{V}_{SE(2)}(\mathbf{x},\theta)$. The circle $S^1$ was sampled with $N_\theta = 32$. All distances were computing via the fast marching algorithm of \cite{Mirebeau2014,mirebeau-hal-01538482} except for the sub-Riemannian approximations, which were computed directly using (\ref{eq:coordinatesfirstkind}) and (\ref{eq:approxSE2}). The position-orientation balancing parameter was set to $\xi=0.01$.

\begin{table}
\caption{Perceptual grouping performance for the 2D retinal image experiments in terms of percentage of correct key point connections (\# of false connections in parenthesis).}
\centering
\begin{tabular}{lll|l}
\toprule
\multicolumn{2}{l}{Distance}                                                      &$C = 1$   & $C \neq 1$     \\
\midrule
$\lVert \mathbf{x} - \mathbf{y} \rVert$             &{\tiny($\mathbb{R}^2\;\;\;\;\;\,|$Eucl.)}  & 89.99\% (362)           & 95.96\%  (146) \\
$d_1(g,h)$                             &{\tiny($SE(2)|$Riem.)}                    & 97.51\% (90)           & 99.64\%  (13) \\
$d_0(g,h)$                             &{\tiny($SE(2)|$Sub-Riem.)\hspace{-0.5em}}                & 99.75\% (9)           & 99.83\%  (6) \\
{\tiny $|\operatorname{Log}(g^{-1}h)|_{\xi,\zeta}$ } \hspace{-1.8em} &{\tiny($SE(2)|\approx$Sub-Riem.)\hspace{-8.0em}}   & 99.72\% (10)    & - \\
\bottomrule
\end{tabular}%
\label{tab:PGPerformanceSE2}
\end{table}

\begin{figure*}
  \centering
  \includegraphics[width=\textwidth]{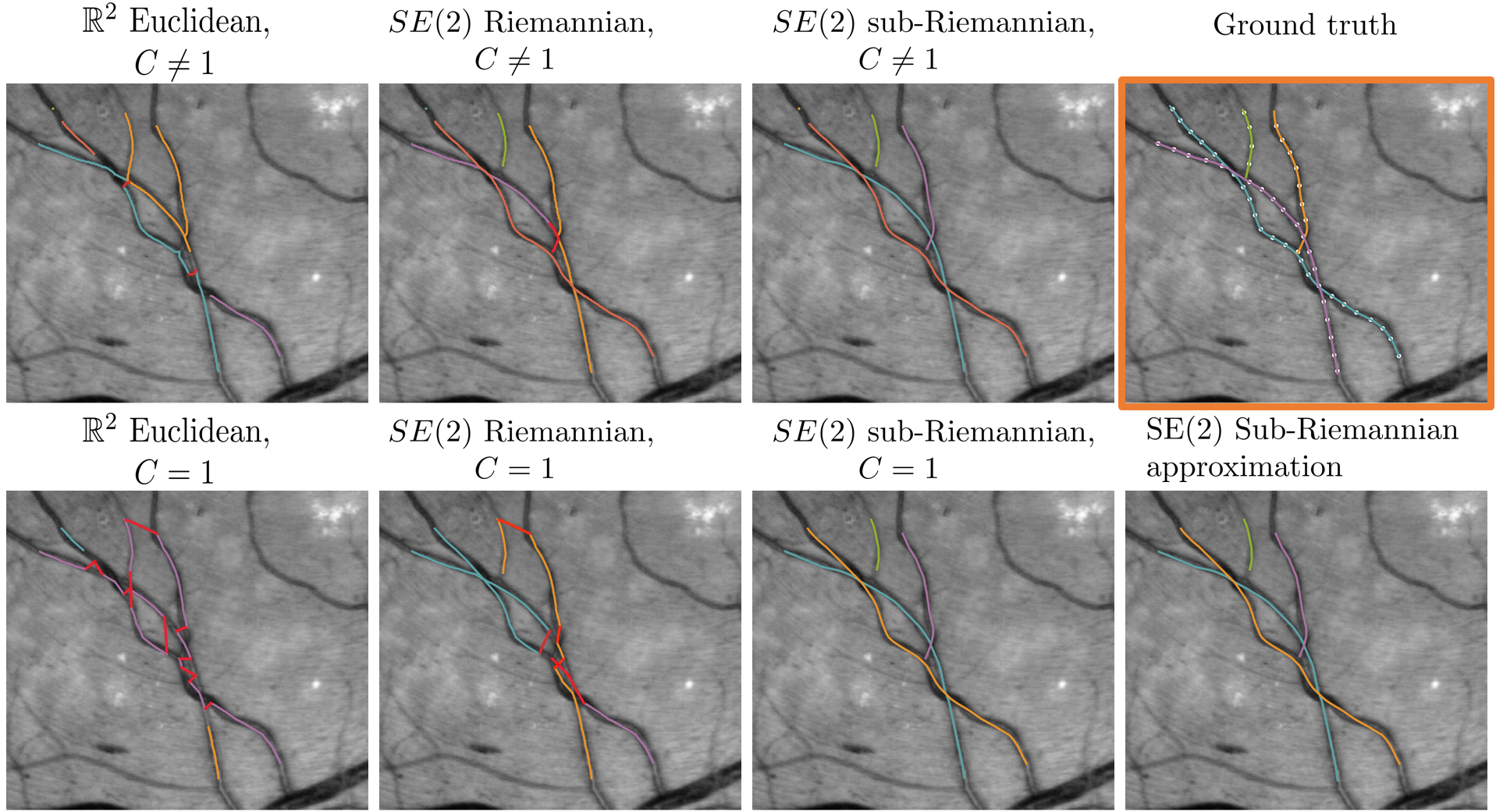}
  \caption{Example 1 of the retinal vessel grouping experiments. Each connected component has its own color (note that the colors might not match between experiments as the number of recovered components may differ), false connections are indicated in red. Top row: experiments with data-adaptive distances ($C\neq1$), and the ground truth vessel components including the automatically generated key points. Bottom row: experiments without data-adaptive distance ($C=1$).}\label{fig:resultsSE2-1}
\end{figure*}

\begin{figure*}
  \centering
  \includegraphics[width=\textwidth]{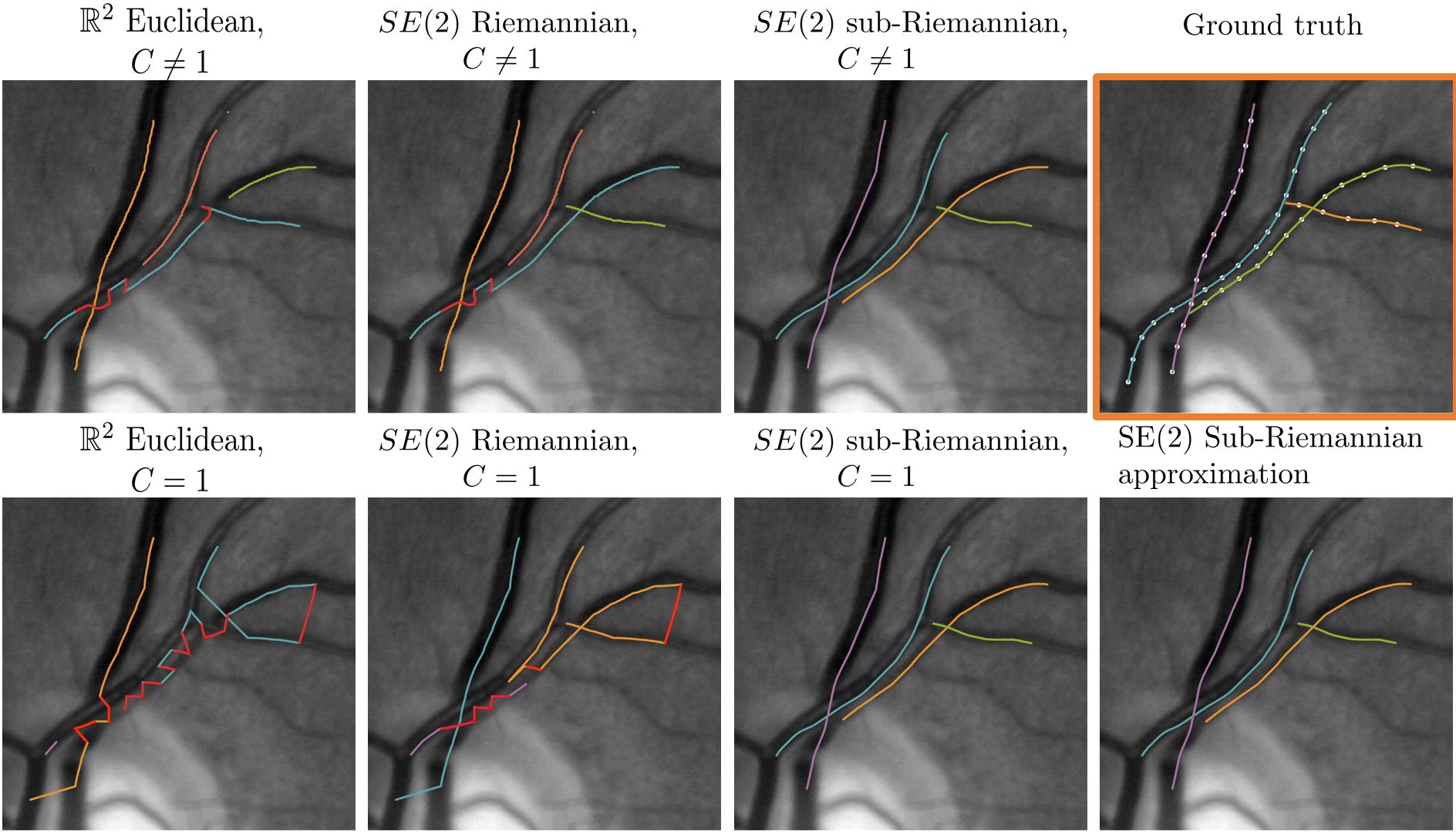}
  \caption{Example 2 of the retinal vessel grouping experiments. Each connected component has its own color (note that the colors might not match between experiments as the number of recovered components may differ), false connections are indicated in red. Top row: experiments with data-adaptive distances ($C\neq1$), and the ground truth vessel components including the automatically generated key points. Bottom row: experiments without data-adaptive distance ($C=1$).}\label{fig:resultsSE2-2}
\end{figure*}

\subsubsection{Results}
Table~\ref{tab:PGPerformanceSE2} gives a quantitative overview of the results, Figs.~\ref{fig:resultsSE2-1} and \ref{fig:resultsSE2-2} shows the results on two of the 52 retina patches. From Table~\ref{tab:PGPerformanceSE2} we make the following observations and conclusions:
\begin{enumerate}
  \item Perceptual grouping is preferred in the lifted domain $SE(2)$ instead of in the based domain $\mathbb{R}^2$. This suggest that taking orientation into account in the grouping is essential.
  \item A sub-Riemannian geometry on $SE(2)$ is preferred over a ($\xi$-isotropic) Riemannian geometry. This suggests that a sub-Riemannian geometry is necessary to deal with the complex geometry at crossings and parallel tracks (cf. Figs.~\ref{fig:resultsSE2-1} and \ref{fig:resultsSE2-2}).
  \item The results obtained with the sub-Riemannian distances on $SE(2)$ for $C=1$ are almost equal. This suggests that the approximations are quite accurate, and that for $C=1$ the analytic approximations may be preferred due to speed and algorithm complexity considerations.
  \item Overall, results for $C\neq 1$ are better than for $C=1$. Note however, that the sub-Riemannian distances on $SE(2)$ for $C=1$ are still better then the Euclidean distance on $\mathbb{R}^2$ and Riemannian distance on $SE(2)$ for $C\neq 1$, and only slightly under performs relative to the sub-Riemannian $C\neq 1$ case. This again shows that sub-Riemannian geometry is preferred, whether data is included in the metric tensors or not.
\end{enumerate}
We conclude that in perceptual grouping of 2D vessels a sub-Riemannian geometry in $SE(2)$ is preferred over a Euclidean geometry in $\mathbb{R}^2$, or a Riemannian geometry in $SE(2)$. When accurate vesselness maps are available, it is preferable to use these in the distances. Furthermore, if one aims to design a easy to implement and efficient perceptual grouping pipeline, approximate sub-Riemannian distances should be used. With only a 2D key point tracking algorithm, a method for estimating orientations, and the analytic approximate distances (\ref{eq:approxSE2}) one obtains very accurate grouping results.

\subsection{Perceptual Grouping in $SE(3)$}\label{subsec:PGSE3}
\subsubsection{Experimental setup}
To quantify and study the influence of different distances in perceptual grouping algorithms for 3D vessels we make use of synthetic 3D images. For these experiments 10 volumes were generated, each with 6 random paths. Each path was generated with a Monte-Carlo simulation of a random walk in $SE(3)$ (see e.g. \cite[Ch.~3.5]{ZhangDuits2014}). Due to the random construction it might occur that 2 paths cross each other. This is physiologically unrealistic (vessels in 3D might bifurcate or touch, but never grow through each other), but it does make the experiments more challenging.

For each volume a binary center line mask was constructed using the generated ground truth paths. The volumes were of size $51\times 51\times 51$ voxels. The distance between key points was set to $l_{max}=5$ voxels. The maximum geodesic arc length distance in the perceptual grouping algorithm was set to $s_{max}=15$ voxels. The orientation at each key point was again estimated as the orientation that gave maximum response in $\mathcal{V}_{SE(n)}$ (Eq.~(\ref{eq:vesselnessSEn})). The sphere $S^2$ was sampled with 200 orientations using Euler angles with $\mathbf{n}(\beta,\gamma) = \mathbf{R}_{\gamma,\beta,\alpha}.\mathbf{e}_z$, with $\beta \in \{\frac{\pi}{2 N_\beta}, 2\frac{\pi}{2 N_\beta}, ...,\pi - \frac{\pi}{2 N_\beta}\}$, $\gamma \in \{0, \frac{\pi}{N_\beta}, ...,2\pi - \frac{\pi}{N_\beta}\}$, with $N_\beta=10$, and with $\mathbf{R}_{\gamma,\beta,\alpha}$ given by (\ref{eq:eulerangles}). In the lifted metric tensor we set $\xi=1$.

\begin{table}
\caption{Perceptual grouping performance for the 3D synthetic volume experiments in terms of percentage of correct key point connections (\# of false connections in parenthesis).}
\centering
\begin{tabular}{lll|l}
\toprule
\multicolumn{2}{l}{Distance}                                                      &$C= 1$   & $C \neq 1$     \\
\midrule
$\lVert \mathbf{x} - \mathbf{y} \rVert$             &{\tiny($\mathbb{R}^3\;\;\;\;\;\,|$Eucl.)}       & 89.99\% (78)           & 97.97\%  (16) \\
$d_1(g,h)$                             &{\tiny($SE(3)|$Riem.)}                    & 93.02\% (54)           & 98.32\%  (13) \\
$d_0(g,h)$                             &{\tiny($SE(3)|$Sub-Riem.)\hspace{-0.5em}}                & 96.79\% (25)           & 98.32\%  (13) \\
{\tiny $|\operatorname{Log}(g^{-1}h)|_{\xi,\zeta}$ } \hspace{-1.8em} &{\tiny($SE(3)|\approx$Sub-Riem.)\hspace{-8.0em}}         &    97.17\%  (22)                 & - \\
\bottomrule
\end{tabular}%
\label{tab:PGPerformanceSE3}
\end{table}

\subsubsection{Results}
\begin{figure*}
  \centering
  \includegraphics[width=1\textwidth]{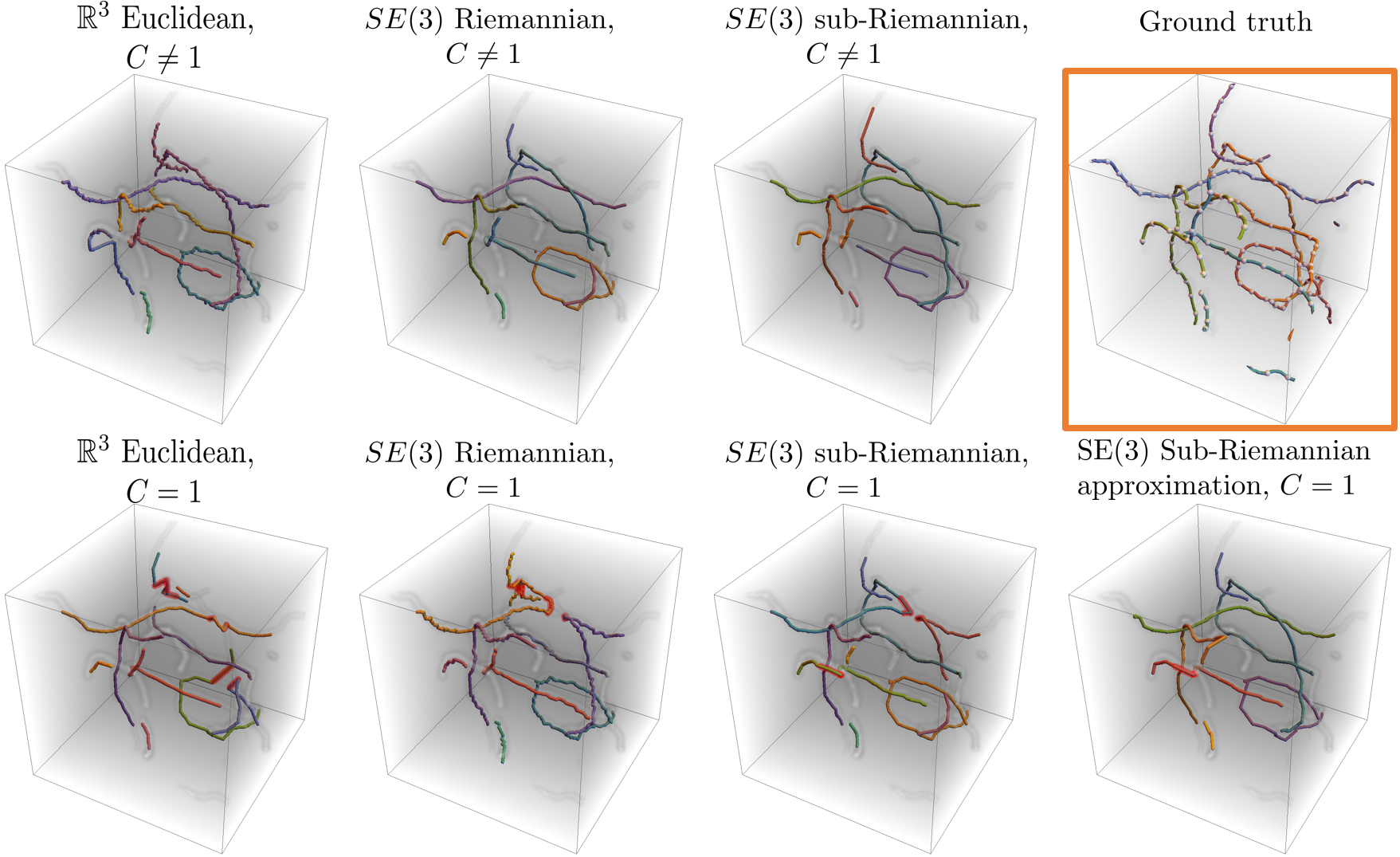}
  \caption{Example 1 of the 3D synthetic vessel grouping experiments. Each connected component has its own color (note that the colors might not match between experiments as the number of recovered components may differ), false connections are indicated in red. Top row: experiments with data-adaptive distances ($C\neq1$), and the ground truth vessel components including the automatically generated key points. Bottom row: experiments without data-adaptive distance ($C=1$).}\label{fig:resultsSE3-1}
\end{figure*}

Table~\ref{tab:PGPerformanceSE3} gives a quantitative overview of the results, Figs.~\ref{fig:resultsSE3-1} shows the results on one of the ten synthetic volumes. From Table~\ref{tab:PGPerformanceSE3} we can draw the same conclusions as for the $SE(2)$ case (using a sub-Riemannian geometry and including data adaptivity improves results). Here, however, we make two additional observations
\begin{enumerate}
  \item Data-adaptive fast marching seems less sensitive to the choice of metric, but tracking in the lifted domain $SE(3)$ still improves results. This can be explained by the fact that the volume is relatively sparse, and by the fact that the cost function $C$ is constructed from ground truth data (the best possible cost). If the cost function dominates the metric, then the intrinsic energy/geometry has a smaller influence.
  \item Out of all $C = 1$ distances (no-data adaptivity) the grouping via the nilpotent distance approximations on $SE(3)$ give best performance, even better then for the true sub-Riemannian distance. This can be explained by the fact that for long distances from the origin, the approximation gradually loose their sub-Riemannian nature and allows more non-horizontal behavior, as in the Riemannian case. It could be that, due to the discrete sampling of the sphere, not all orientations are accurately estimated. The grouping based on the sub-Riemannian distance approximations seems less sensitive to such errors.
\end{enumerate}

\section{Conclusion}\label{sec:conclusion}
In this paper we have proposed an efficient approach for perceptual grouping of local orientations via nilpotent approximations of sub-Riemannian distances in the roto-translation group $SE(n)$. The quantitative experiments on grouping of retinal blood vessels in 2D images, and perceptual grouping in challenging 3D synthetic volumes, showed that 1) sub-Riemannian geometry is essential in achieving top performance and 2) that the grouping approach via the fast analytic approximations performs almost equally, or better, than the data-adaptive fast marching approaches.

The sub-Riemannian distances on $SE(2)$ and $SE(3)$ were approximated via norms on exponential coordinates of the first kind (obtained via the logarithmic map). In both quantitative and visual comparison it was found that the approximations accurately follow the true sub-Riemannian distances, a conclusion which was further supported by the equal performance in quantitative perceptual grouping experiments. We also numerically showed that the
weighted logarithmic norms used in this paper provide a more accurate approach for approximating the heat kernel and fundamental solution of the sub-Laplacian on $SE(n)$, compared to previous approaches \cite{Duits2010,portegies2016arxiv,citti2006cortical}.

Since the sub-Riemannian distance approximations are analytic, they are easy to implement and fast to compute. An interesting line of further research would be to embed the sub-Riemannian distances in other algorithms that rely on the quantification of the distance between local orientations. The results of this paper could be further improved by augmenting the sub-Riemannian distances with additional features (like cross-sectional profile descriptors) and use a global graph optimization approach as in \cite{turetken2016reconstructing,estrada2015retinal}. The potential of using sub-Riemannian distances in such problems is demonstrated by the experiments of this paper.

\begin{acknowledgements}
The following are gratefully acknowledge for their influence on the manuscript: Remco Duits at Eindhoven University of Technology for suggestions on logarithmic approximations of the heat kernel and the fundamental solution on $SE(2)$, $SE(3)$ and the Heisenberg group; Laurent Cohen at University Paris Dauphine for fruitful discussions on geodesic methods and perceptual grouping; Jean-Marie Mirebeau at Laboratoire de math{\'e}matiques d’Orsay, Universit{\'e} Paris-Saclay for providing efficient and generic fast marching code.
The reviewers are gratefully acknowledged for their valuable feedback on the manuscript.  The research leading to the results of this article has received funding from the European Research Council under the European Communitys 7th Framework Programme (FP7/20072014)/ERC grant agreement No. 335555 (Lie Analysis).
\end{acknowledgements}

\appendix
\section{Optimization of the Folland-Kaplan-Kor{\'a}nyi gauge parameter $\zeta$}\label{app:optimizationZeta}
In Figs.~\ref{fig:DistancesSE2} and Fig.~\ref{fig:DistancesSE3} we visually compared the nilpotent approximations of the sub-Riemannian distance on $SE(2)$ and $SE(3)$ respectively. In this appendix we support by means of a quantitative comparison our choices for $\zeta=44$ and $\zeta=100$ which appear in the logarithmic approximations (Folland-Kaplan-Kor{\'a}nyi gauge) of $|\operatorname{Log}(g^{-1}h)|_{\xi,\zeta}$ as defined in (\ref{eq:approxDistSE2}) and (\ref{eq:approxDistanceSE3}) on respectively $SE(2)$ and $SE(3)$.

\subsection{Optimization of $\zeta$ for the $SE(2)$ approximations}
\begin{figure}
  \centering
  \includegraphics[width=0.48\textwidth]{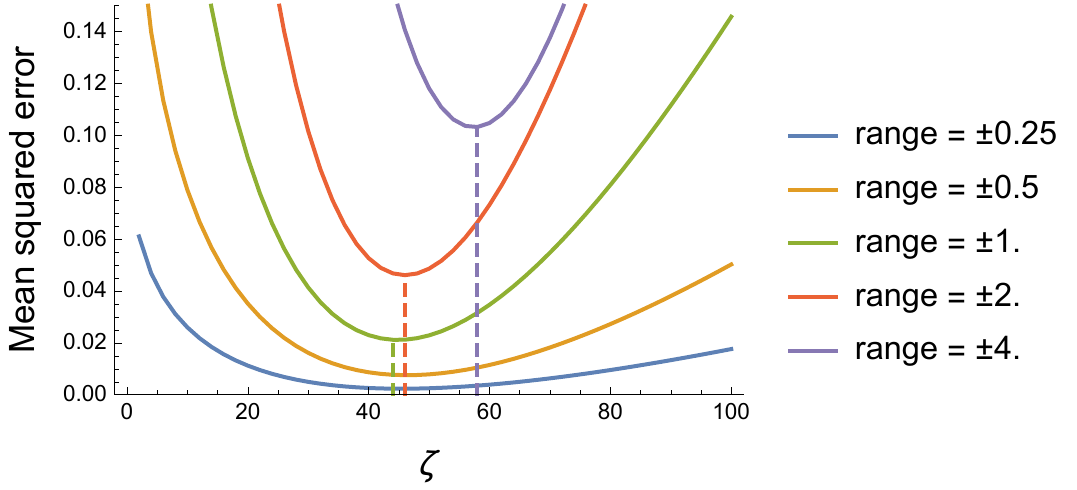}
  \caption{{\color{red}
   Mean squared errors between the sub-Riemannian distance on $SE(2)$ (see Eq.~(\ref{eq:SRDistancesSE2})) and its approximation (see Eq.~(\ref{eq:approxDistSE2})). The errors is computed for varying choices of $\zeta$ and on a varying grid size (from $x,y\in [-0.5,0.5]$ to $x,y\in [-4,4]$).}}\label{fig:NumericalCompSE2}
\end{figure}

In the quantitative comparison on $SE(2)$ we computed the $\mathbb{L}_2$ error between $d_0(g,h)$ and the approximation $|\operatorname{Log}(g^{-1}h)|_{\xi,\zeta}$ with $\xi=1$ on a grid with a varying spatial domain size (from $x,y\in [-0.5,0.5]$ to $x,y\in [-4,4]$), and with varying choices of $\zeta$. The results are shown in Fig.~\ref{fig:NumericalCompSE2} and are computed as follows.

The reference sub-Riemannian distance $d_0$ on $SE(2)$ was computed once via an anisotropic fast marching algorithm \cite{Mirebeau2014,mirebeau-hal-01538482} on a grid which sampled $x,y \in [-4,4]$ at a sub-pixel resolution of $0.01$ with 128 orientations. The numerically computed sub-Riemannian distance volume was thus of dimensions $801 \times 801 \times 128$.

In each experiment (with fixed spatial range and $\zeta$) the squared error between $d_0$ and its approximation was sampled on a regular grid that covered the specified domain with $41\times 41 \times 128$ points. The averaged errors are plotted in Fig.~\ref{fig:NumericalCompSE2}.

Here we see that the approximation becomes more accurate towards the origin ($x,y\in [-0.5,0.5]$) and that parameter $\zeta$ has to be chosen larger in order to keep the anisotropy for longer distances from the origin. The choice $\zeta=44$ generally gave the best approximations and we rely on this setting in the experiments on $SE(2)$.

\subsection{Optimization of $\zeta$ for the $SE(3)$ approximations}
\begin{figure}
  \centering
  \includegraphics[width=0.48\textwidth]{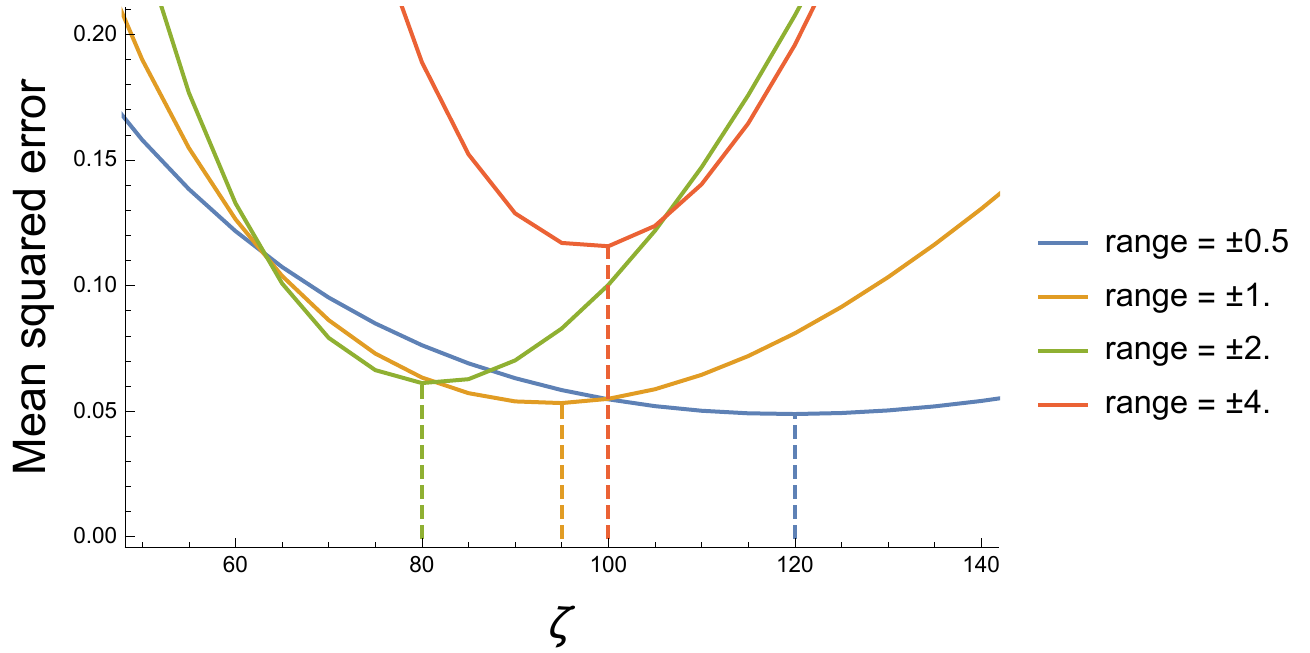}
  \caption{{\color{red}
   Mean squared errors between the sub-Riemannian distance on $SE(3)$ (see Eq.~(\ref{eq:SRDistancesSE2}) and (\ref{eq:SRMetricTensorSE3})) and its approximation (see Eq.~(\ref{eq:approxDistanceSE3})). The errors is computed for varying choices of $\zeta$ and on a varying grid size (from $x,y,z\in [-0.5,0.5]$ to $x,y,z\in [-4,4]$).}}\label{fig:NumericalCompSE3}
\end{figure}

In the quantitative comparison on $SE(3)$ we computed the $\mathbb{L}_2$ error between $d_0(g,h)$ and the approximation $|\operatorname{Log}(g^{-1}h)|_{\xi,\zeta}$ with $\xi=1$ on a grid with a varying spatial domain size (from $x,y,z\in [-0.5,0.5]$ to $x,y,z\in [-4,4]$), and with varying choices of $\zeta$. The results are shown in Fig.~\ref{fig:NumericalCompSE3} and are computed as follows.

The reference sub-Riemannian distance $d_0$ on $SE(3)$ was also computed once via an anisotropic fast marching algorithm \cite{Mirebeau2014,mirebeau-hal-01538482} on a grid which sampled $x,y,z \in [-4,4]$ at a sub-pixel resolution of $0.1$ with $31\times62$ Euler angles (cf. Sec.~\ref{subsec:PGSE3}). The numerically computed sub-Riemannian distance volume was thus of dimensions $201 \times 201 \times 201 \times 31 \times 62$.

In each experiment (with fixed spatial range and $\zeta$) the squared error between $d_0$ and its approximation was sampled on a regular grid that covered the specified domain with $21\times 21 \times 21 \times 31 \times 62$ points. The averaged errors are plotted in Fig.~\ref{fig:NumericalCompSE3}.

Here we see that the approximation becomes more accurate towards the origin ($x,y,z\in [-0.5,0.5]$). However, compared to the $SE(2)$ experiments we do see a less stable localization of the optimal parameter $\zeta$ with varying spatial resolutions. This behaviour can be explaind by 1) the sub-Riemannian distances are numerically computed via a fast marching algorithm using Euler angles (which do not uniformly sample the sphere), and 2) the spatial resolution of the computed reference sub-Riemannian distance volume was only $0.1$ (due to computer memory constraints). Although very accurate from an application point of view, the sub-Riemannian distances on $SE(3)$ are not exact, and the numerical errors induced by the algorithm may explain the variation in optimal $\zeta$ (in particular for the region close to the origin). Overall, the choice $\zeta=100$ seems to be reasonable in all ranges, and this was confirmed by visual comparison of the distance maps in Fig.~\ref{fig:DistancesSE3}.

\bibliographystyle{spmpsci}      
\bibliography{references}

%
%

\end{document}